\documentclass{article}
\usepackage{amssymb}
\usepackage{amsmath}
\usepackage{enumerate}
\usepackage{amscd}
\usepackage{graphicx}

\usepackage{latexsym}
\usepackage{epsfig}
\usepackage[latin1]{inputenc}

\def\div{\mathfrak{Div}}
\def\deg{\mathfrak{Deg}}

\def\r{\mathbb{R}}
\def\n{\mathbb{N}}
\def\c{\mathbb{C}}

\def\s{\mathbb{S}}
\def\d{\mathbb{D}}

\def\z{\mathbb{Z}}

\def\k{\mathbb{K}}

\def\nb{\mathcal{N}}

\def\Fg{\mathfrak{F}}

\def\Rg{\mathfrak{R}}
\def\Wg{\mathfrak{W}}

\def\Sg{\mathfrak{S}}

\newenvironment{proof}{\trivlist
\item[\hskip\labelsep{\em Proof}\,:]}{\hfill{$\Box$}\endtrivlist}

\title{\huge Uniform Approximation by Complete Minimal Surfaces of Finite Total Curvature in $\r^3$}
\author{\Large Francisco J. L\'{o}pez   \thanks{Research partially supported by MCYT-FEDER research projects MTM2007-61775 and MTM2011-22547, and Junta de Andaluc\'{i}a Grant P09-FQM-5088.
Complete minimal surfaces of finite total curvature, compact minimal surfaces, Runge's Theorem, Mergelyan's Theorem.}}

\newtheorem{lemma}{Lemma}[section]
\newtheorem{remark}{Remark}[section]
\newtheorem{theorem}{Theorem}[section]
\newtheorem{proposition}{Proposition}[section]\newtheorem{corollary}{Corollary}[section]
\newtheorem{assertion}{Claim}[section]
\newtheorem{definition}{Definition}[section]
\begin{document}
\maketitle

\begin{abstract}We prove that any compact minimal surface in $\r^3$ can be uniformly approximated by complete minimal surfaces of finite total curvature in $\r^3$ is obtained. This Mergelyan's type result can be extended to the family of complete minimal surfaces  of weak finite total curvature, that is to say, having  finite total curvature on  regions of finite conformal type. We deal only with the orientable case.
\end{abstract}
\section{Introduction} \label{sec:intro}

The classical theorems of Mergelyan and Runge deal with the uniform approximation problem for holomorphic functions on planar regions  by rational functions on the complex plane.  They  extend to   interpolation problems and approximation results of continuous functions on Jordan curves and meromorphic functions on regions, among other applications. Specially interesting is the approximation  by  meromorphic functions with prescribed  zeros and poles on compact Riemann surfaces ({\em algebraic} approximation). For instance, see  the works by Bishop \cite{bishop}, Scheinberg \cite{sche1,sche2} and Royden \cite{roy} for a good setting.

These results have  played an interesting role in the general theory of minimal surfaces, taking part in very sophisticated arguments for constructing complete (or proper) minimal surfaces  that are far from being algebraic in any sense (see the pioneering works by Jorge-Xavier \cite{jorge-xavier}, Nadirashvili \cite{nadi} and Morales \cite{morales}).

Complete minimal surfaces of finite total curvature ({\em FTC} for short)  have occupied a relevant position in the global theory of minimal surfaces  since its origin. Progress in this area has depended, in a essential manner, on their special analytic and geometric properties.  Huber \cite{huber} proved that
 if $M$ is a Riemann surface with possibly non empty compact boundary $\partial (M)$ which admits a  conformal complete minimal immersion with FTC  in $\r^3,$ then $M$ has {\em finite conformal type} accordingly with the following
 \begin{definition} \label{def:finconf}
A Riemann surface $M$ with possibly $\partial (M) \neq \emptyset$ is said to be of finite conformal type if $M$ is conformally equivalent to $M^c-E,$ where $M^c$ is a compact Riemann surface and $E\subset M^c-\partial(M^c)=M^c-\partial (M)$ is a  finite set (the  topological ends of $M$).
\end{definition}  
In addition, R. Osserman  \cite{osserman} showed that the Weierstrass data of  such a immersion extend meromorphically to the so-called Osserman compactification $M^c$ of $M.$   It is also worth  mentioning that any finitely punctured compact Riemann surface admits a  (non-rigid) conformal complete minimal immersion of FTC in $\r^3,$ see Pirola \cite{pirola}.
 
This paper is devoted to developing a general approximation theory for minimal surfaces in $\r^3.$ In this context, complete minimal surfaces of FTC will play the same role as rational functions in complex analysis. Among other things, we prove that any compact minimal surface can be uniformly approximated by complete minimal surfaces of FTC   (see Theorem II below). This cousin result  of Runge's and Mergelyan's theorems  establishes a natural connection between the local and global theories of minimal surfaces,  and leads to  natural geometrical applications  (bridge constructions, isoperimetric inequalities,  immersing problems, and general existence theorems for minimal surfaces).  Moreover, our  methods allow  control over the conformal structure and flux map of the approximating surfaces.  Further developments can be found, for instance, in \cite{a-l,a-f-l,exotic}. 

For a thorough exposition of these results, the following notations are required.

A conformal complete minimal immersion $X:M \to \r^3$ is said to be of {\em weak finite total curvature} ({\em WFTC} for short) if $X|_\Omega$ has FTC for all  regions $\Omega \subset M$ of finite conformal type.

If $X:M \to \r^3$ is a conformal minimal immersion and $\gamma\subset M$ is an oriented closed curve,  the flux  of $X$ on $\gamma$ is given by $p_X(\gamma):=\int_\gamma \mu(s) ds,$ where   $s$ is an oriented arc length parameter on $\gamma$  and $\mu(s)$ the corresponding conormal vector of $X$ at $\gamma(s)$ for all $s.$ Recall that $\mu(s)$ is the unique unit  tangent vector of $X$ at $\gamma(s)$ such that $\{d X(\gamma'(s)),\mu(s)\}$ is a positive basis. Since $X$ is a harmonic map, $p_X(\gamma)$ depends only on the homology class of $\gamma$ and the well defined flux map $p_X:\mathcal{H}_1(M,\z) \to \r^3$ is a group homomorphism.

As usual, a topological surface is said to be {\em open} if it is non-compact and has empty topological boundary. In the sequel, $\nb$ will denote an arbitrary but fixed open Riemann surface.

\begin{definition}\label{def:basica} Let $M$ be a 
proper subset of $\nb$ all whose connected components are regions of $\nb.$ We denote by $\mathcal{M}(M)$  the space of conformal complete minimal immersions $X:M \to\r^3$ of WFTC, extending as a conformal minimal immersion to a neighborhood of $M$ in $\nb.$
\end{definition}
If $M$ consists of finitely many pairwise disjoint regions of finite conformal type,  $ \mathcal{M}(M)$ is the space of  conformal complete minimal immersions of $M$ in  $\r^3$ with  FTC extending beyond $M$ in $\nb.$  The space $ \mathcal{M}(M)$ will be endowed with the topology of the {\em uniform convergence on (not necessarily compact) regions of finite conformal type in $M.$}

Our main result is the following theorem (see Theorems \ref{th:density} and \ref{th:parabo}):

\begin{quote} 
{\bf Theorem I (Fundamental Approximation Theorem): }{\em Let $M$ be a finite collection of pairwise disjoint regions in $\nb$ of finite conformal type, and let $\beta$ be a finite collection of compact analytical Jordan arcs (possibly some of them closed Jordan curves) in $\nb$ meeting each other finitely many times (the cases $M=\emptyset$ or $\beta=\emptyset$ are allowed). Assume that $\overline{\beta-M}$ has finitely many connected components, set $S=M\cup \beta,$ and assume that $\nb-S$ contains no relatively compact connected components.

Then, for any  two smooth conformal maps\footnote{See Definition \ref{def:smooth}.} $X:S \to \r^3,$ $N:S\to \s^2,$ and group homomorphism   $q:\mathcal{H}_1(\nb,\z) \to \r^3$  satisfying:
\begin{itemize}
\item $X|_M\in \mathcal{M}(M)$ and $X|_\beta$ is an immersion,
\item $N|_M$ is the Gauss map of $X|_M$ and $N|_\beta$ is  normal to $X|_\beta,$ and
\item $q|_{\mathcal{H}_1(S,\z) }=p_{X},$
\end{itemize}
 there  exists $\{Y_n\}_{n \in \n} \subset  \mathcal{M}(\nb)$  such that $\{(Y_n|_{S}-X,{N_n}|_{S}-N)\}_{n \in \n} \to (0,0)$ uniformly on $S,$ where $N_n$ is the Gauss map of $Y_n,$ and $p_{Y_n}=q$ for all $n.$}
\end{quote}
 As a corollary, $ \mathcal{M}(\nb)\neq \emptyset$  for any open Riemann surface $\nb,$ and even more, for any group homomorphism $q:\mathcal{H}_1(\nb,\z)\to \r^3$  we can find $Y \in  \mathcal{M}(\nb)$ with $p_Y=q.$ Choosing $\nb$ of finite conformal type and  $q=0,$ one obtains Pirola's theorem \cite{pirola} as a corollary.

The Fundamental Approximation Theorem can also be used in  general connected sum constructions for complete minimal surfaces of FTC (see \cite{exotic}).  Other results of this kind can be found in Kapouleas \cite{kapou} and  Yang \cite{yang}. 
Perhaps, the most basic and useful consequence of Theorem I is the following corollary, in which all the involved immersions are of FTC:
\begin{quote} 
{\bf Theorem II (Basic Approximation Theorem): }{\em Assume that $\nb$ has finite conformal type, and   let $M\subset \nb$ be a compact region such that $\nb-M$ contains no relatively compact components .
 
Then, for any  $X\in  \mathcal{M}(M)$ there exists  $\{Y_n\}_{n \in \n}\subset  \mathcal{M}(\nb)$ such that $\{Y_n|_M-X\}_{n \in \n} \to 0$ uniformly on $M$ and $p_{Y_n}|_{\mathcal{H}_1(M,\z)}=p_X$ for all $n.$}
\end{quote}

The paper is laid out as follows. Section \ref{sec:pre} is devoted to some preliminary results on Algebraic Geometry, and in Section \ref{sec:wei-spin} we go over the Weierstrass and spinorial representations of minimal surfaces. Section \ref{sec:appro} contains the main results: the Fundamental Approximation Theorem is proved in Subsection \ref{subsec:algbrid} in the particular case when $\nb$ has finite conformal type and $S$ is {\em isotopic} to $\nb,$ whereas its general version and the Basic Approximation Theorem  are obtained in Subsection \ref{subsec:parabo}.

In a forthcoming paper \cite{lop} the author will extend this analysis to the non orientable case.

\section{Preliminaries on Riemann surfaces} \label{sec:pre}

As usual, we call $\overline{\c}=\c \cup \{\infty\}$ the extended complex plane or Riemann sphere.  

Let $M$ be a topological surface with possible non-empty topological boundary. As usual, we write  $\partial(M)$ for the one dimensional topological manifold determined by the boundary points of $M,$ and $\mbox{Int}(M)$ the open surface  $M-\partial(M).$ A subset $\Omega\subset M$ is said to be {\em proper} if the inclusion map $j:\Omega\to M$ is a proper topological map, which  simply means that $\Omega$ is closed in $M.$ A proper connected subset  $\Omega\subset M$ is said to be a {\em region} if, endowed with the induced topology, it is a  topological surface with possibly non-empty boundary.  An open connected subset of $\mbox{Int}(M)$ will be called a {\em domain} of $M.$   Given $S \subset M,$ we write $S^\circ$ and $\overline{S}$ for the topological interior and  closure of $S$  in $M.$

\begin{definition}\label{def:annular}
Given two regions $\Omega$ and $\Omega^*$ of $M$ with finitely many boundary components,   $\Omega^*$ is said to be an  {\em extension} of $\Omega$ if $\Omega$ is a proper subset of $\Omega^*,$ $\Omega\cap \partial(\Omega^*)=\emptyset$ and $\Omega^*-\Omega^\circ$ contains no compact connected components disjoint from $\partial(\Omega^*).$ In particular the induced group homomorphism $j_*:\mathcal{H}_1(\Omega,\r) \to \mathcal{H}_1(\Omega^*,\r)$ is injective, where  $j:\Omega \to \Omega^*$ is the inclusion map (up to the natural identification we consider $\mathcal{H}_1(\Omega,\r) \subset \mathcal{H}_1(\Omega^*,\r)$).

Likewise,  $\Omega^*$ is said to be an   {\em annular extension} of $\Omega$ if $\Omega^*$  is an extension of $\Omega$ and $\Omega^*-\Omega^\circ$ consists of finitely many compact annuli and once punctured closed discs. If in addition  $\Omega^*-\Omega^\circ$ is compact (that is to say, it is a finite collection of compact annuli), then $\Omega^*$ is said to be a {\em trivial annular extension or closed tubular neighborhood} of $\Omega$ in $M.$ In this case $\Omega$ and $\Omega^*$ are homeomorphic.
\end{definition}
These notions can be extended to the case when $\Omega$ and $\Omega^*$ are a finite collection of pairwise disjoint regions in $\nb.$\\

Assume now that $M$ is a Riemann surface. 

For any $W\subset M,$ we denote by $\div(W)$  the free commutative group of divisors of $W$ with multiplicative notation. If $D=\prod_{i=1}^n Q_i^{n_i} \in \div(W),$ where $n_i \in \z-\{0\}$ for all $i,$ the set $\{Q_1,\ldots,Q_n\}$ is said to be the {\em support} of $D.$ We denote by $\deg:\div(W) \to \z$  the  group homomorphism given by the degree map $\deg(\prod_{j=1}^t Q_j^{n_j})=\sum_{j=1}^t n_j.$ A divisor $\prod_{i=1}^n Q_i^{n_i} \in \div(W)$ is said to be {\em integral} if $n_i\geq 0$ for all $i.$ Given $D_1,$ $D_2 \in \div(W),$ $D_1 \geq D_2$ if and only if $D_1 D_2^{-1}$ is  integral.

Let $W$ be an open subset of $\mbox{Int}(M),$ and let $F:W\to\c$ be a meromorphic function. We denote by $(F)_0$ and $(F)_\infty$ the integral divisors of zeros and poles of $F$ in $W,$  respectively, and call $(F)=(F)_0/(F)_\infty$ the divisor of $f$ in $W.$ If $V\subset W$ is  a subset (normally,  a region or a finite collection of them) and $f=F|_V,$  we also write $(f)_0$ and $(f)_\infty$ for the corresponding integral divisors of zeros and poles of $F$ in $V,$  respectively, and call $(f)=(f)_0/(f)_\infty$ the divisor of $f$ in $V.$ Likewise for meromorphic 1-forms.

\subsection{Compact Riemann surfaces}
The background of the following results can be found, for instance, in \cite{farkas}. 

In the sequel, $R$ will denote a compact  Riemann surface with genus $\nu \geq 1$ and empty boundary. We denote by $\Wg_m(R)$ and $\Wg_h(R)$ the spaces of meromorphic and holomorphic 1-forms on $R,$ respectively, and call $\Fg_m(R)$ the space of meromorphic functions on $R.$

 Label $\mathcal{H}_1(R,\z)$ as the $1^{st}$  homology group with integer coefficients of $R.$ Let $B=\{a_j,b_j\}_{j=1,\ldots,\nu}$ be a canonical homology basis of $\mathcal{H}_1(R,\z),$ and write $\{\xi_j\}_{j=1,\ldots,\nu}$ the associated dual basis of  $\Wg_h(R),$   that is to say, the one satisfying that $\int_{a_k} \xi_j=\delta_{jk}, \quad j,\;k=1,\ldots, \nu.$ 
 
 Denote by  $\Pi=(\pi_{jk})_{j,\,k=1,\ldots,\nu}$ the Jacobi period matrix with entries $\pi_{jk}=\int_{b_k} \xi_j, \quad j,\;k=1,\ldots, \nu.$ This matrix is symmetric and has positive definite imaginary part. We denote by $L(R)$ the lattice over $\z$ generated by the $2\nu$-columns of the $\nu \times 2 \nu$ matrix $(I_\nu,\Pi),$ where $I_\nu$ is the identity matrix of dimension $\nu.$

If  $f \in \Fg_m(R)$ and $(f)_0$ and $(f)_\infty\in \div(R)$ are the  integral divisors of zeros and poles of $f$ in $R,$ respectively, we call $(f)=(f)_0/(f)_\infty$  the principal divisor associated to $f.$ Likewise, if $(\theta)_0$ and $(\theta)_\infty\in \div(R)$ are the integral divisors of zeros and poles of  $\theta \in \Wg_m(R),$ respectively,  we call $(\theta)=(\theta)_0/(\theta)_\infty$ the canonical divisor of $\theta.$

 Finally, set $J(R)=\c^\nu/L(R)$ the Jacobian variety of $R,$ which is a compact, commutative,  complex, $\nu$-dimensional Lie group. 
 Fix $P_0 \in R,$ denote by $\varphi_{P_0}: \div(R) \to J(R), \quad \varphi_{P_0}(\prod_{j=1}^s Q_j^{n_j})= \sum_{j=1}^s n_j \,{}^t(\int_{P_0}^{Q_j} \xi_1,\ldots,\int_{P_0}^{Q_j} \xi_\nu)$ the Abel-Jacobi map with base point $P_0,$ where ${}^t(\, \cdot \,)$ means matrix transpose. If there is no room for ambiguity, we simply write $\varphi.$ 
 
Abel's theorem asserts that  $D\in \div(R)$ is the principal divisor associated to a meromorphic function $f \in \Fg_m(R)$ if and only if $\deg(D)=0$ and $\varphi(D)=0.$  Jacobi's theorem says that $\varphi:R_\nu \to J(R)$ is surjective and has maximal rank (hence a local biholomorphism) almost everywhere, where $R_\nu$ denotes the space of integral divisors in $\div(R)$ of degree $\nu.$   

Riemann-Roch theorem says that $r(D^{-1})=\deg(D)-g+1+i(D) $ for any $D \in \div(R),$  where $r(D^{-1})$ (respectively, $i(D)$) is the dimension of the complex vectorial space  of functions  $f\in \Fg_m(R)$ (respectively, 1-forms $\theta \in \Wg_m(R)$) satisfying that $(f) \geq D^{-1}$ (respectively, $(\theta) \geq D$).

By Abel's theorem, the point $\kappa_R:=\varphi((\theta)) \in J(R)$ does not depend on $\theta\in \Wg_m(R).$  It is called  the {\em  vector of the Riemann constants}. 
Write ${\cal S}(R)$ for the set containing the $2^{2 \nu}$ solutions of the algebraic equation $2 s=\kappa_R$ in $J(R).$ Any element of ${\cal S}(R)$ is said to be a {\em spinor structure} on $R.$  A  1-form $\theta\in \Wg_m(R)$ is said to be {\em spinorial} if $(\theta)=D^2$ for a divisor $D \in \div(R).$  Denote by $\Sg_m(R)$ (respectively, $\Sg_h(R)$)  the set of spinorial meromorphic (respectively, spinorial holomorphic)  1-forms on $R.$ Two  1-forms $\theta_1,$ $\theta_2\in \Sg_m(R)$ are said to be spinorially equivalent, written $\theta_1 \sim \theta_2,$ if there exists $f \in \Fg_m(R)$ such that $\theta_2=f^2 \theta_1.$  Notice that a class $\Theta \in \frac{\Sg_m(R)}{\sim}$  determines a unique spinor structure $s_\Theta\in {\cal S}(R).$ Indeed, it suffices to take $\theta \in \Theta$ and define $s_\Theta=\varphi(D),$ where $D\in \div(R)$ is determined by the equation  $D^2=(\theta).$ By Abel's theorem $s_\Theta$ does not depend on the chosen $\theta \in \Theta.$   

The map $\frac{\Sg_m(R)}{\sim}\rightarrow {\cal S}(R),$ $\Theta \mapsto s_\Theta$ is bijective. To see this, take $s\in {\cal S}(R)$ and use Jacobi's theorem to find an integral divisor $D' \in \div(R)$ of degree $\nu$ satisfying $\varphi(D')=s.$ By Abel's theorem, $(D' P_0^{-1})^2$ is the canonical divisor associated to a spinorial meromorphic 1-form whose corresponding class $\Theta_s$ in $\frac{\Sg_m(R)}{\sim}$  satisfies $s_{\Theta_s}=s$ (as indicated above, $P_0$ is the initial condition of the Abel-Jacobi map). 

Spinor structures can be also introduced in a more topological way. Indeed, take  $s\in {\cal S}(R)$ and  $\theta \in \Theta_s.$ For any embedded loop $\gamma \subset R,$ consider an open annular neighborhood $A$ of $\gamma$ and a conformal parameter $z:A \to \{z \in \c\;:\; 1<|z|<r\}.$ Set  $\xi_s(\gamma)=0$ if $\sqrt{\theta(z)/dz}$  has  a well defined branch on $A$ and  $\xi_s(\gamma)=1$ otherwise, and note that this number does not depend on the chosen annular conformal chart.   The induced map $\xi_s:\mathcal{H}_1(R,\z) \to \z_2$ does not depend on $\theta \in \Theta_s$ and defines a group homomorphism. Furthermore, $\xi_{s_1}=\xi_{s_2}$ if and only if $s_1=s_2,$ and therefore ${\cal S}(R)$ can be identified with the set of group morphisms  $\mbox{Hom}(\mathcal{H}_1(R,\z),\z_2).$ We simply write $\xi_\Theta=\xi_{s_\Theta},$ for any $\Theta \in \frac{\Sg_m(R)}{\sim}.$\\

\subsection{Riemann surfaces of finite conformal type}
Let $M$ be a Riemann surface of  finite conformal type with possibly $\partial (M)\neq \emptyset$ (see Definition \ref{def:finconf}), and write $M=M^c-\{E_1,\ldots,E_a\},$ where $M^c$ is compact and 
$\{E_1,\ldots, E_a\}\subset  M^c-\partial (M^c).$ The compact Riemann surface $M^c$ is said to be the {\em Osserman compactification} of $M$ (uniquely determined up to biholomorphisms). 
Any compact Riemann surface is of finite conformal type (in this case, the set of topological ends is empty).

Attaching a conformal disc to each connected component of $\partial (M^c)=\partial (M),$ we get a compact Riemann surface $R$ without boundary that will be called  a  {\em conformal compactification} of $M.$ With this language, $M^c=R-(\cup_{j=1}^b U_j),$ where $U_1,\ldots,U_b$ are open discs in $R$ with pairwise disjoint closures. Notice that  $R$ depends on the gluing process of the conformal discs, hence conformal compactifications of $M$  are not unique. 

As usual, call $\mathcal{H}_1(M,\z)$ the $1^{st}$ homology group of $M$ with integer coefficients.  

Set  $\Sg_m(M)$ the space of meromorphic 1-forms $\theta$ on $\mbox{Int}(M)$  satisfying that
\begin{itemize}
\item any zero or pole of $\theta$ in $\mbox{Int}(M)$ has even order, and
\item $\theta$ extend meromorphically to $\mbox{Int}(M^c).$
\end{itemize}
In a similar way, we call  $\Sg_h(M)$ the space of  $\theta \in \Sg_m(M)$ such that $\theta$ is holomorphic on $\mbox{Int}(M).$

Two 1-forms $\theta_1,$ $\theta_2 \in \Sg_m(M)$ are said to be {\em spinorially equivalent} if there exists a meromorphic function  $f$ on $\mbox{Int}(M^c)$ such that $\theta_2=f^2 \theta_1.$ As above, we define the map  
$$ \xi: \frac{\Sg_m(M)}{\sim} \rightarrow \mbox{Hom}(\mathcal{H}_1(M,\z),\z_2),\quad  \Theta \mapsto \xi_\Theta.$$

\begin{lemma} \label{lem:spin}
The map $\xi: \frac{\Sg_m(M)}{\sim} \rightarrow \mbox{Hom}(\mathcal{H}_1(M,\z),\z_2)$ is bijective. 
\end{lemma}
\begin{proof} Standard monodromy arguments show that $\xi$ is injective. 

Put $\partial (M)=\cup_{j=1}^b c_j,$ where $c_j$ is a Jordan curve for all $j$ and $c_{j_1} \cap c_{j_2}=\emptyset$ when $j_1 \neq j_2.$ Consider a family $V_1,\ldots,V_a$ of pairwise disjoint closed discs in $M^c-\partial(M)$ such that $E_i \in {V_i}^\circ$ for all $i=1,\ldots,a.$  Label $r:=a+b>0$ and $\{d_1,\ldots,d_r\}=\{\partial(V_i),\;i=1,\ldots,a\} \cup \{c_j,\; j=1,\ldots,b\}.$ 
Let $R$ be a conformal compatification of $M,$ and fix a homology basis $\{a_1,\ldots,a_\nu,b_1,\ldots,b_\nu\}$  of $\mathcal{H}_1(R,\z).$ 

We know that $\{a_1,\ldots,a_\nu,b_1,\ldots,b_\nu,d_1,\ldots,d_{r-1}\}$ is a basis of $\mathcal{H}_1(M,\z),$   so $\mbox{Hom}(\mathcal{H}_1(M,\z),\z_2)$ contains $2^{2 \nu+r-1}$ elements.  
Write  $\frac{\Sg_m(R)}{\sim}=\{\Theta_j,\: j=1,\ldots 2^{2 \nu}\}.$ Choose $\theta_j \in \Theta_j$ for each $j$ and call $f_j=\theta_j/\theta_1 \in \Fg_m(R),$ $j=1,\ldots,  2^{2 \nu}.$ Since $\Theta_i$ and $\Theta_j$ correspond to different spinor structures on $R,$ $i \neq j,$  $\sqrt{\theta_i/\theta_j}$ has no well defined branches on $R,$ hence the same holds on $\mbox{Int}(M).$ Thus $\{\theta_j|_{\mbox{Int}(M)}\;:\; j=1,\ldots 2^{2 \nu}\}$ are pairwise spinorially inequivalent in $\Sg_m(M).$ Write $M^c=R-\cup_{j=1}^b U_j,$ where $U_j$ is an open disc in $R$ with $\partial (U_j)=c_j$ for all $j,$ and fix $E_{a+j} \in {U_j},$ $j=1,\ldots,b.$ For any $J \subseteq \{1,\ldots,r-1\},$ $J \neq \emptyset,$ use Jacobi's theorem to find  an integral divisor $D_J \in \div(R)$ of degree $\nu$ verifying $$\varphi(D_J^2 P_0^{2\sharp(J)-2}E_r^{-\sharp(J)}\prod_{j \in J}E_j^{-1} )=\kappa_R,$$  where $\sharp(J)$ is the cardinal of $J$ and $P_0$ is the initial condition   of $\varphi.$ By Abel's theorem, there exists $\tau_J\in \Wg(R)$ with canonical divisor $(\tau_J)= D_J^2 P_0^{2\sharp(J)-2} E_r^{-\sharp(J)}\prod_{j \in J}E_j^{-1}.$  Since $f_i \tau_J/\theta_j$ has a pole of  odd order at some $E_h,$ $h\in \{1,\ldots,r\},$ $(f_i \tau_J)|_{\mbox{Int}(M)}$ and $\theta_j|_{\mbox{Int}(M)}$ are not spinorially equivalent in $\Sg_m(M),$ $i,j \in \{1,\ldots,2^{2 \nu}\},$ and likewise for any pair $(f_{i_1} \tau_{J_1})|_{\mbox{Int}(M)},$ $(f_{i_2} \tau_{J_2})|_{\mbox{Int}(M)}$  with $(i_1,J_1) \neq (i_2,J_2).$ Thus   $\{\theta_j|_{\mbox{Int}(M)}, \;j=1,\ldots,2^{2 \nu}\} \cup \{(f_i \tau_J)|_{\mbox{Int}(M)},\; i=1,\ldots, 2^{2 \nu},\; J \subseteq \{1,\ldots,r-1\},\; J \neq \emptyset\}$ contains 
$2^{2 \nu+r-1}$ pairwise spinorially inequivalent 1-forms in $\Sg_m(M),$ proving that $\xi$ is surjective.

\end{proof}

\subsection{Approximation results on Riemann surfaces} \label{sec:appprelim}

In this section we recall some basic approximation theorems in complex analysis. 

We first adopt some conventions and fix some notations.

\begin{remark}
In the sequel,  $\nb$  will denote an open Riemann surface.
\end{remark}

\begin{definition} \label{def:ends}
We denote by $\nb^c$ the Riemann surface obtained by filling out all the conformal punctures of $\nb$ (that is to say, the annular ends of $\nb$ of finite conformal type). In other words, $\nb^c$ is the union of the  Osserman compactifications of all regions in $\nb$ of  finite conformal type. 
\end{definition}

Given  $V\subset \nb,$ a connected component $U$ of $\nb-V$ is said to be {\em bounded} if $\overline{U}$ is compact.

If $V\subset \nb$ is an arbitrary  subset, we denote by $V^c$ the subset of $\nb^c$ obtained by attaching to $V$ the isolated points of  $\nb^c-V$ (that is to say, the conformal punctures of $V$).

Let us introduce the special subsets of $\nb$ on which our later constructions are based.

\begin{figure}[ht]
\begin{center}
\includegraphics[width=7cm,height=3.5cm]{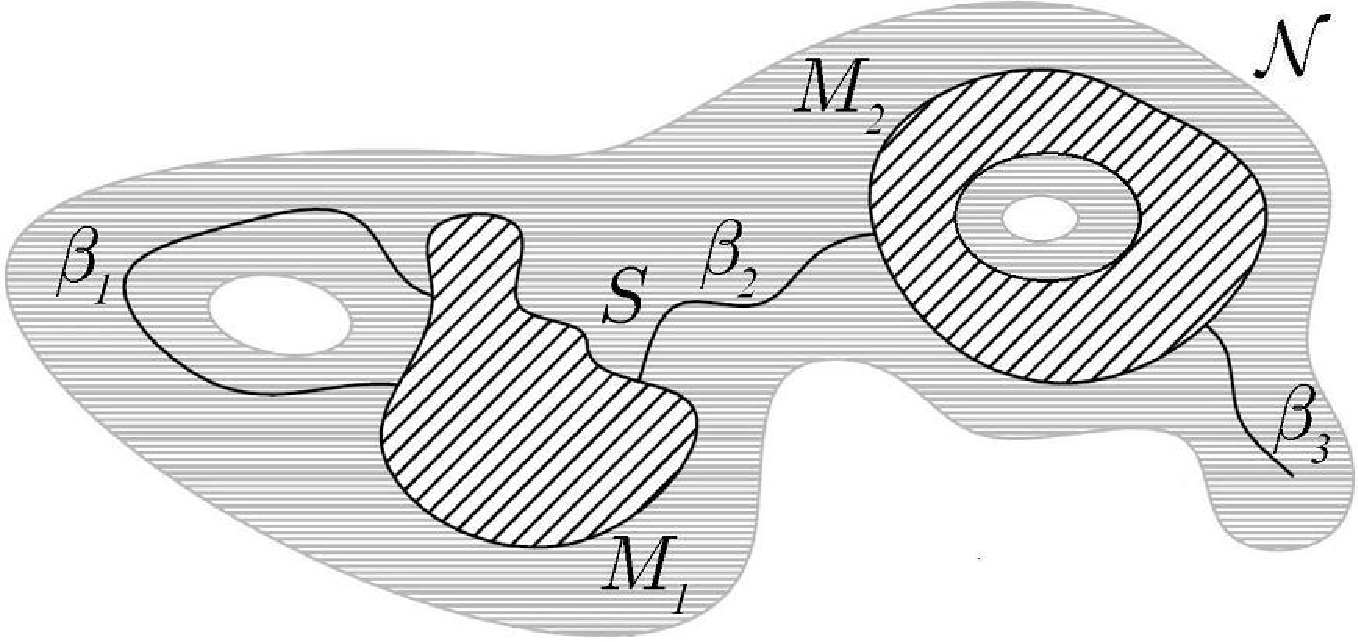} \caption{An admissible subset $S=M\cup \beta$ with $k=2$ and $m=3.$}  \label{fig:dibu1} 
\end{center}
\end{figure}

\begin{definition} \label{def:admissible}
A proper subset $S\subset \nb,$ $S\neq \emptyset,$ is said to be {\em admissible} in $\nb$ if it admits a decomposition $S=M\cup \beta,$ where 
\begin{itemize}
\item  $M$ is either empty or consists of finitely many pairwise disjoint regions  $M_1,\ldots,$ $M_k,$ $k\geq 1,$  of finite conformal type and non-empty boundary, 
\item  $\beta$ is either empty or consists of finitely many analytical compact Jordan arcs $\beta_1,\ldots, \beta_m$ in $\nb,$  possibly some of them closed Jordan curves,
\item  $\{\beta_i\cap\beta_j\,|\; i\neq j\}$ is finite and $\overline{\beta-M}$ consists of finitely many compact Jordan arcs (possibly some of them closed Jordan curves), and 
\item $\nb-S$ has no bounded components. 
\end{itemize}
If $S$ is admissible in $\nb,$ we call 
\begin{itemize}
\item $\partial(S):=\partial (M) \cup \beta,$ and 
\item $S^c=M^c\cup \beta\subset \nb^c$ (the Osserman compactification of $S$).   
\end{itemize}
 
See Figure \ref{fig:dibu1}.
\end{definition}

\begin{remark} \label{re:beta0}
When $\beta_j$ is not a closed Jordan curve, we always suppose that $\beta_j\subset \beta_{0,j},$ where $\beta_{0,j}$ is either an open analytical arc or a closed curve. We make the convention $\beta_{0,j}=\beta_j$ if  $\beta_j$ is a closed Jordan curve. Furthermore, we will assume that $\{\beta_{0,i}\cap\beta_{0,j}\,|\; i\neq j\}$ is finite and $\overline{\beta_0-M}$ has finitely many connected components as well, where $\beta_0=\cup_{j=1}^m \beta_{0,j}.$ 
\end{remark} 
Notice that if $S$ is admissible in $\nb,$ then  $S^c\cup \nb$ is an {\em open Riemann surface} and $S^c$ is admissible in $S^c\cup \nb.$

\begin{figure}[ht]
\begin{center}
\includegraphics[width=8.5cm,height=3.2cm]{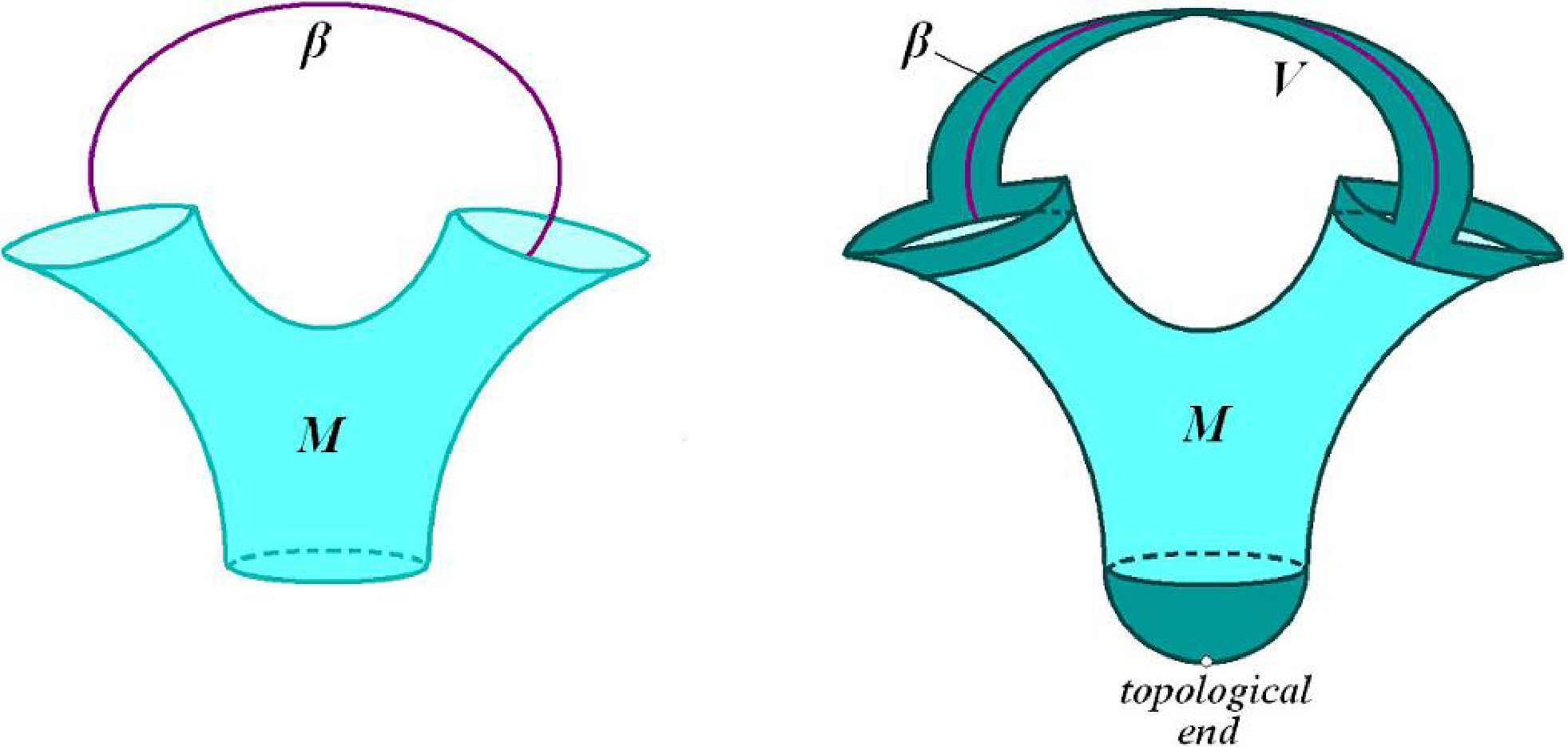} \caption{An annular extension $V$ of $M\cup \beta$ when $k=m=1.$}  
\label{fig:dibu2} 
\end{center}
\end{figure}

A region $V \subset \nb$ is said to be an {\em annular extension} of an admissible $S=M\cup\beta$ in $\nb$  if  it is a annular extension of a small closed tubular neighborhood $S_0$ of $S$ in $\nb$ (which can be defined in the standard way with the help of a  complete Riemannian metric on $\nb$). In particular, $S\subset V^\circ,$ any relatively compact connected component of  $V-(M\cup \beta)$ meets $ \partial(V),$  $V-(M\cup \beta)$ consists of a finite collection of conformal annulus and conformal  once punctured discs, and the induced homomorphism  $j_*:\mathcal{H}_1(S,\z) \to \mathcal{H}_1(V,\z)$ is an isomorphism, where $j:S \to V$ is the inclusion map. See Figure \ref{fig:dibu2} and  Definition \ref{def:annular}. If in addition the closure of  $V$ is a closed tubular  neighborhood of $S_0$ (that is to say, $V-(M\cup \beta)$ contains no conformal once punctured discs), then $V$ is said to be a {\em closed tubular  neighborhood}  of $S$ ($S_0$ itself is a closed tubular  neighborhood of $S$).

Throughout this paper, we will only deal with analytical objects extending meromorphically to conformal punctures. With this philosophy in mind, we need the following spaces of functions and 1-forms.

\begin{definition}
Let $V$ be a finite collection of pairwise disjoint regions or domains in $\nb.$ We denote by
\begin{itemize}
\item $\Fg_m(V)$  the space of meromorphic  functions on $V$ extending meromorphically to $V^c$ (if $V$ is a region, we always suppose that $f$ extends meromorphically  to a neighborhood of $V$ in $\nb$),
\item $\Fg_h(V)$  the space of  holomorphic functions on $V$ extending meromorphically to $V^c$ (if $V$ is a region, we always suppose that $f$ extends  holomorphically to a neighborhood of $V$ in $\nb$), and
\item $\Fg_h(V^c)$ the space of functions in $\Fg_h(V)$ extending holomorphically to $V^c.$ 
\end{itemize}

Likewise, we call  $\Wg_m(V),$ $\Wg_h(V),$ and $\Wg_h(V^c)$ the analogous spaces of 1-forms. 
\end{definition}
For instance, $\Wg_h(\nb)$ is the space of holomorphic  1-forms on $\nb$ extending meromorphically to $\nb^c.$
The inclusions $\Fg_h(V^c)\subset \Fg_h(V)\subset \Fg_m(V)$ and $\Wg_h(V^c)\subset \Wg_h(V)\subset \Wg_m(V)$ are trivial.\\ 

\begin{remark}
The space $\Fg_m(V^c)$ of functions in $\Fg_m(V)$ extending meromorphically to $V^c$ is nothing but $\Fg_m(V),$ so in most cases we will avoid this redundant notation. Likewise for $\Wg_m(V^c).$ 
\end{remark}

Let us present now the corresponding spaces of functions and 1-forms on admissible sets in $\nb$ and their natural topologies. 
In the remaining section, $S=M\cup \beta \subset \nb$ will be an admissible subset in $\nb.$

\begin{definition} \label{def:spacfun}
We call $\Fg_m(S)$ (respectively, $\Fg_h(S)$) the space of continuous  functions $f:S \to \overline{\c}$ such that $f|_M\in \Fg_m(M)$  (respectively, $f|_M\in \Fg_h(M)$) and
$f(P) \neq \infty$ for all $P \in \beta.$

The space of functions $f\in \Fg_h(S)$ extending holomorphically to $M^c$ will be labeled by $\Fg_h(S^c).$

\end{definition}

In a natural way, $\Fg_h(S^c) \subset \Fg_h(S)\subset \Fg_m(S).$ These spaces are endowed with the topology of the uniform convergence on $S^c$ (or equivalently, on $S$), also called the $\mathcal{C}^0(S)$-topology.  

\begin{definition}
We shall say that a function $f \in \Fg_m(S)$  can be uniformly approximated on $S$ by functions in $\Fg_m(\nb)$ if there exists  a sequence $\{f_n\}_{n \in \n}$ in $\Fg_m(\nb)$    such that
$\{f_n|_S\}_{n\in \n}\to f$ in the $\mathcal{C}^0(S)$-topology, that is to say,  $\{|f_n|_S-f|\}_{n \in \n} \to 0$ uniformly on $S^c.$ In this case, $f_n-f\in\Fg_h(S^c)$ for all $n\in \n,$ and in particular, all $f_n,$ $n\in \n,$ have the same set of poles as $f$  on $S^c.$ 

If in addition $f\in \Fg_h(S)$ and  $\{f_n\}_{n \in \n}\subset \Fg_h(\nb)$(respectively,  $f\in \Fg_h(S^c)$ and  $\{f_n\}_{n \in \n}\subset \Fg_h(S^c\cup \nb)$), one gets the corresponding notion of uniform approximation by holomorphic functions in $\Fg_h(\nb)$ (respectively,  in $\Fg_h(S^c\cup \nb)$). 
\end{definition}

A complex $1$-form $\theta$ on  $S$ is said to be of type $(1,0)$ if for any conformal chart $(U,z)$ on $\nb,$ one has $\theta|_{U \cap S}=f(z) dz$ for some  $f:U \cap S \to \overline{\c}.$

\begin{remark}
Fix an auxiliary complete conformal Riemannian metric $\rho_\nb^2$ on $\nb^c.$
\end{remark}

\begin{definition}
We call $\Wg_m(S)$ (respectively, $\Wg_h(S)$) the space of  1-forms $\theta$ of type $(1,0)$  such that 
$\theta|_M\in \Wg_m(M)$  (respectively, $\theta|_M\in \Wg_h(M)$), $\theta(P) \neq \infty$ for all $P \in \beta,$ and  $\theta|_\beta$ is continuous.

The space of 1-forms  $\theta \in \Wg_h(S)$ extending holomorphically to $M^c$ will be labeled by $\Wg_h(S^c).$
\end{definition}

In a natural way, $\Wg_h(S^c)\subset \Wg_h(S)\subset \Wg_m(S).$ These spaces are endowed with the topology of the uniform convergence on $S^c$ (or equivalently, on $S$),  also called the $\mathcal{C}^0(S)$-topology.  
The convergence $\{\theta_n\}_{n\in \n} \to \theta$ must be understood as $\{|\frac{\theta_n-\theta}{\rho_\nb}|\}_{n\in\n}\to 0$ uniformly on $S^c.$

\begin{definition}
We shall say that a 1-form $\theta$ in $\Wg_m(S)$ can be uniformly approximated on $S$ by 1-forms in $\Wg_m(\nb)$ if there exists a sequence $\{\theta_n\}_{n \in \n}$ in $\Wg_m(\nb)$  such that $\{|\frac{\theta_n-\theta}{\rho_\nb}|\}_{n \in \n} \to 0$ in the $\mathcal{C}^0(S)$-topology.  In this case, $\theta_n-\theta\in\Wg_h(S^c)$ for all $n\in \n,$ and in particular, all $\theta_n,$ $n\in \n,$ have the same set  of poles as $\theta$  on $S^c.$

If in addition $\theta\in \Wg_h(S)$ and  $\{\theta_n\}_{n \in \n}\subset \Wg_h(\nb)$(respectively,  $\theta\in \Wg_h(S^c)$ and  $\{\theta_n\}_{n \in \n}\subset \Wg_h(S^c\cup \nb)$), one gets the corresponding notion of uniform approximation by holomorphic 1-forms in $\Wg_h(\nb)$ (respectively,  in $\Wg_h(S^c\cup \nb)$). 
\end{definition}
Notice that these notions of convergence for 1-forms do not depend on the auxiliary conformal metric $\rho_\nb^2$ in $\nb^c.$

E. Bishop \cite{bishop}, H. L. Royden \cite{roy} and S. Scheinberg \cite{sche1,sche2}, among others, have proved several extensions of  Runge's and Mergelyan's theorems. For our purposes, we need only the following compilation result:

\begin{theorem} \label{th:runge} Let $S=M\cup\beta \subset \nb$ be an admissible subset in $\nb.$ Then  any  $f \in \Fg_m(S)$ can be uniformly approximated on $S$  by functions $\{f_n\}_{n \in \n}$ in $\Fg_m(\nb) \cap {\Fg_h}(\nb-{\cal P}_f),$  where ${\cal P}_f=f^{-1}(\infty)\subset M.$ Furthermore, if $D \in \div(M)$ is an integral divisor,  then the approximating sequence$\{f_n\}_{n\in \n}$ can be chosen  so that $\big(f|_M-f_n|_{M}\big)_0 \geq D.$
\end{theorem}
 
\begin{remark} \label{re:runge}
In most applications of Theorem \ref{th:runge}, the divisor $D$ is chosen satisfying that $D\geq (f|_M)_0.$ 
If $f$ never vanishes on $\partial(S)$ and $D\geq (f|_M)_0,$ then $\frac{f_n|_S}{f} \in \Fg_h(S^c)$ for all $n$ and  $\{\frac{f_n|_S}{f}\}_{n\in \n} \to 1$ in the ${\cal C}^0(S)$-topology.  
\end{remark}

\section{Analytic Representations of Minimal Surfaces}\label{sec:wei-spin}

Let us review some basic facts about minimal surfaces. 

Fix an open Riemann surface $\nb$ and an auxiliary complete conformal Riemannian metric $\rho_\nb^2$ on $\nb^c,$ and keep the notations of Section \ref{sec:appprelim}.

Let $M$ denote a finite collection of pairwise disjoint regions in $\nb$. Endow $\mathcal{M}(M)$ (see Definition \ref{def:basica})  with the following ${\cal C}^0(M)$-topology:

\begin{definition} \label{def:convergencia}
A sequence $\{X_n\}_{n \in \n} \subset  \mathcal{M}(M)$ is said to converge in the ${\cal C}^0(M)$-topology to  $X_0 \in \mathcal{M}(M)$ if for any region $\Omega \subset M$ of finite conformal type, $\{X_n|_\Omega-X_0|_\Omega\}_{n \in \n} \to 0$ uniformly on $\Omega,$ that is to say, in the topology associated to the norm of the supremum on $\Omega.$ In particular,  $X_n-X_0$ extends harmonically to $\Omega^c$  by Riemann's removable singularity theorem for all $n,$ and $\{(X_n-X_0)|_{\Omega^c}\}_{n \in \n} \to 0$  in the norm of the maximum on $\Omega^c.$

If $M$ has finite conformal type, this topology coincides with the one of the uniform convergence on $M.$
\end{definition}

Let $X=(X_j)_{j=1,2,3}$ be a conformal minimal immersion in $\mathcal{M}(M).$ Write $\partial_z X_j=\phi_j$  and notice that $\partial_z X_j\in \Wg_h(M)$ for all $j.$ Since $X$ is conformal and minimal,  then  $\phi_1=\frac{1}{2}(1/g-g) \phi_3$ and $\phi_2=\frac{i}{2}(1/g+g) \phi_3,$ where $g \in \Fg_m(M)$ is, up to the stereographic projection, the Gauss map of $X.$ The pair $(g,\phi_3)$ is known as the {\em Weierstrass representation} of $X.$

Clearly $X(P)=X(Q)+\mbox{Re} \int_{Q}^P (\phi_1,\phi_2,\phi_3),$ $P,$ $Q \in M.$ The induced intrinsic metric $ds^2$ on $M$ and its Gauss  curvature $\mathcal{K}$ are given by the expressions: 
\begin{equation} \label{eq:regular}
ds^2=\sum_{j=1}^3 |\phi_j|^3=\frac{1}{4} |\phi_3|^2(\frac{1}{|g|}+|g|)^2,\quad \mathcal{K}=-\left(\frac{4 |dg||g|}{ |\phi_3| (1+|g|^2)^2}\right)^2.
\end{equation}
The total curvature of $X$ is given by $c(X):=\int_M \mathcal{K} dA,$ where  $dA$ is the area element of $ds^2,$ and the {\em flux map} of $X$ by  the expression $p_X:\mathcal{H}_1(M,\z) \to \r^3,$ $p_X(\gamma)=\mbox{Im} \int_\gamma \partial_z X.$

If $\partial(M)$ is compact and  $X$ is a minimal complete immersion of FTC, Huber, Osserman and Jorge-Meeks results \cite{huber,osserman,jorge-meeks} imply that $X$ is proper,  $M$ has finite conformal type, the Weierstrass data $(g,\phi_3)$ of $X$ extend meromorphically to $M^c,$ and the vectorial 1-form $\partial_z X$ has poles of order $\geq 2$ at the ends (i.e., the points of $M^c-M$). 
 
\begin{remark} If $\{X_n,\; n \in \n\} \cup \{X\} \subset \mathcal{M}(M)$ and  $\{X_n\}_{n \in \n}\to X$ in the $\mathcal{C}^0(M)$-topology,  then  the Weierstrass data of $X_n$ converge uniformly to the ones of $X$ on compact regions of $M^c.$ Indeed, just observe that $\{X_n-X\}_{n \in \n}\to 0,$ and so $\{\partial_z X_n-\partial_z X\}_{n \in \n}\to 0,$ uniformly on compact regions of $M^c.$
\end{remark}

Assume now that $M\subset \nb$ is a region of finite conformal type, consider $X \in \mathcal{M}(M),$  and write $(g,\phi_3)$ for its Weierstrass data. Since $ds^2$ has no singularities on $M$ (see equation (\ref{eq:regular})),  $\eta_1=\frac{\phi_3}{g}$ and $\eta_2=\phi_3 g\in  \Sg_h(M^*),$ are spinorially equivalent  in  $\Sg_h(M^*),$ and have no common zeros on $M^*,$ where $M^*$ is any closed tubular  neighborhood of $M$ in $\nb$ to which $X$ extends.   Furthermore, we know that at least one of them has a pole (of order $\geq 2$) at  each puncture in $M^c-M.$  The next lemma shows that the converse is true:

\begin{lemma}[Spinorial Representation] \label{lem:spinrepre} 

Let $M$ be a region in $\nb$ of finite conformal type, and let $M^*$ be a closed tubular  neighborhood of $M$ in $\nb.$ Let  $\eta_1,$ $\eta_2$ be two spinorially equivalent 1-forms in $\Sg_h(M^*)$ such that $|\eta_1|+|\eta_2|$ never vanishes in $M,$  at least one of the 1-forms $\eta_j,$ $j=1,2,$ has a pole at  each point of $M^c-M,$ and $\frac{1}{2} (\eta_1-\eta_2),$ $\frac{i}{2}(\eta_1+\eta_2)$ and $\sqrt{\eta_1 \eta_2}$ have no real periods on $M.$ 

Then the map  $X:M \to \r^3,$
\begin{equation} \label{eq:spinor}
X(P)=\mbox{Re} \int_{P_0}^P (\phi_1,\phi_2,\phi_3),\quad P_0 \in M,
\end{equation}
where $(\phi_j)_{j=1,2,3}=\left(\frac{1}{2} (\eta_1-\eta_2), \frac{i}{2}(\eta_1+\eta_2),\sqrt{\eta_1 \eta_2}\right),$ is well defined and lies in $\mathcal{M}(M).$ 
\end{lemma}
\begin{proof}
Since $\eta_1$ and $\eta_2$ are spinorially equivalent in $\Sg_h(M^*)$ (and obviously lie in $\Wg_m(M)$),  there is $g \in \Fg_m(M)$ such that $\eta_2=g^2\eta_1,$ and therefore $\phi_3:=\sqrt{\eta_1 \eta_2}$ is well defined.  As $\frac{1}{2} (\eta_1-\eta_2),$ $\frac{i}{2}(\eta_1+\eta_2)$ and $\phi_3$ have no real periods on $M,$ then  $X$ is well defined. Furthermore, from our hypothesis $\frac{1}{4} |\phi_3|^2(\frac{1}{|g|}+|g|)^2$ never vanishes on $M,$ hence $X$ is the minimal immersion with  Weierstrass data $(g,\phi_3).$ Following Osserman \cite{osserman},  $X$ is complete and of FTC. 
\end{proof}

The pair $(\eta_1,\eta_2)$ will be called as the {\em spinorial representation} of $X$ (see \cite{kusner} for a good setting).

\subsection{Minimal surfaces on admissible subsets}

We are going to  introduce the natural notion of conformal minimal immersion on an admissible subset of $\nb$ into $\r^3.$ These surfaces will be the initial  conditions for our main problem, that is to say, the natural objects to which we will later approximate by conformal minimal immersions of WFTC on $\nb.$

\begin{remark} 
In the sequel,  $S=M\cup\beta$ will be an admissible subset in $\nb.$ We  use the notations of Definition \ref{def:admissible}, and accordingly to Remark \ref{re:beta0},  consider an analytical extension $\beta_{0,j}$ of   $\beta_j$ for all  $j=1,\ldots,m$.
\end{remark}

\begin{definition} \label{def:smooth}
A map $X:S \to \k,$ where $\k=\c^n,$ $\r^n$ or $\s^n,$ $n\in \n,$ is said to be a {\em smooth conformal map}  if 
\begin{itemize}
\item there exist an open neighborhood $M_0$ of $M$ in $\nb$ and a smooth conformal map $X_0:M_0\to \k$ such that  $X_0|_{M}=X|_M,$
\item there exists a smooth map  $X_j:\beta_{0,j}\to \k$  such that $X_0|_{\beta_{0,j}\cap M_0}= X_j|_{\beta_{0,j}\cap M_0}$ and $X_j|_{\beta_j}=X|_{\beta_j}$  for all $j,$
\item for any intersection point $P\in \beta_j\cap \beta_i,$ $j\neq i,$ either $d(X|_{\beta_{j}})_P= d(X|_{\beta_{i}})_P=0$ or  $d(X|_{\beta_{j}})_P$ and $d(X|_{\beta_{i}})_P\neq0,$ and in the last case  $$\text{$\frac{(\rho_\nb^2)_P(v_i,v_i)}{\|d(X|_{\beta_{i}})_P(v_i)\|^2}=\frac{(\rho_\nb^2)_P(v_j,v_j)}{\|d(X|_{\beta_{j}})_P(v_j)\|^2}$ and $\measuredangle_\nb(v_j,v_i)=\measuredangle(d(X|_{\beta_{j}})_P(v_j),d(X|_{\beta_{i}})_P(v_i)),$}$$ 
where  $v_j$ and $v_i$ are any tangent vectors at $P$ of $\beta_j$ and $\beta_i,$ $\|\cdot\|$ and $\measuredangle$ are the norm and the oriented angle with respect to the Euclidean metric in $\k,$ and  $\measuredangle_\nb$ is the oriented angle in the Riemannian surface $(\nb,\rho_\nb^2).$ 
\end{itemize}
Notice that this notion does not depend on the chosen conformal metric $\rho_\nb^2$ on $\nb,$ and observe that if $P\in \beta_i\cap\beta_j\cap \beta_h$ then $d(X|_{\beta_{i}})_P,$ $d(X|_{\beta_{j}})_P,$ and $d(X|_{\beta_{h}})_P\}$ lie in a plane of the real tangent space at $X(P)$ of $\k.$
\end{definition} 

\begin{definition}
We denote by ${\cal M}(S)$ the space of smooth conformal maps $X:S \to \r^3$ such that $X_j:=X|_{M}\in \mathcal{M}(M)$   and $X|_{\beta}$ is a regular map (or an immersion). It is clear that $Y|_{S}\in {\cal M}(S)$  for all $Y \in \mathcal{M}(\nb).$ 
\end{definition}

The Gauss map has played a fundamental role for the understanding of the conformal geometry of minimal surfaces.  For this reason, it is natural to {\em mark} the immersions $X \in \mathcal{M}(S)$ with  a normal field along $\beta.$  

\begin{figure}[ht]
\begin{center}
\includegraphics[width=6cm,height=4.5cm]{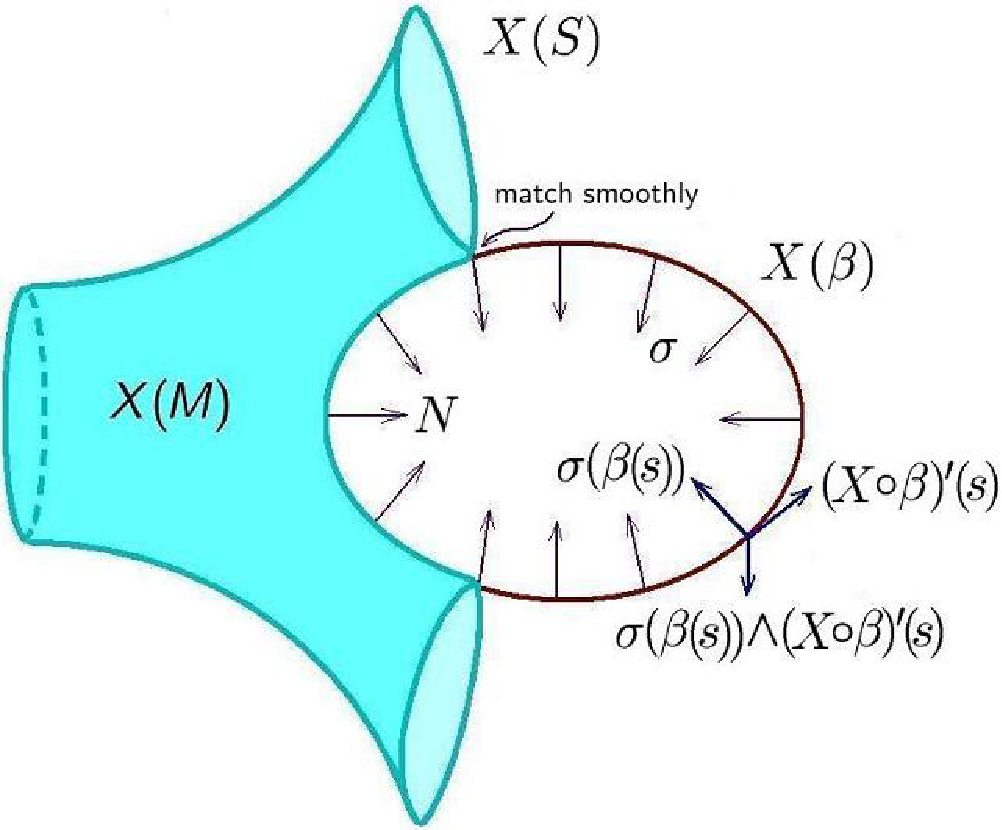} \caption{A smooth normal field $\sigma$ with respect to $X$  along $\beta.$ }  \label{fig:dibu4} 
\end{center}
\end{figure}

\begin{definition} Take $X \in {\cal M}(S),$  and let $N:M\to \s^2$ denote the Gauss map of $X|_M.$
 A  map $\sigma:\beta \to \s^2$ is said to be {\em a smooth normal field} with  respect to $X$ along $\beta$ if  $\sigma({\beta_{j}}(t))$ is orthogonal to $(X\circ \beta_{j})'(t)$ for any smooth parameter $t$ on $\beta_{j}$ and for all $j,$ and the map $$\text{$N_\sigma:S \to \s^2,$ $N_\sigma|_M=N,$ $N_\sigma|_\beta=\sigma,$}$$ is smooth and conformal accordingly to Definition \ref{def:smooth}.

By definition, $N_\sigma$ is said to be the {\em generalized Gauss map} of the marked immersion $(X,\sigma).$
\end{definition}
See Figure \ref{fig:dibu4}.

The following space of immersions will be crucial.
\begin{definition}
We call  $\mathcal{M}^*(S)$  as the  space of marked immersions  $X_\sigma:=(X,\sigma),$ where $X \in {\cal M}(S)$ and $\sigma$ is a  smooth normal field with respect to $X$ along $\beta.$
For any  $X_\sigma,$ $Y_\varpi \in \mathcal{M}^*(S),$ set $$\|X _\sigma-Y_\varpi\|_{1,S}=\|X-Y\|_{0,S}+\|N_\sigma-N_\varpi\|_{0,S},$$  where  $\|\cdot\|_{0,S}$ means $\sup_S \|\cdot\|$ and $\|\cdot\|$ is the  Euclidean norm in $\r^3.$

We endow  $\mathcal{M}^*(S)$ with the $\mathcal{C}^1(S)$-topology of the uniform convergence  of  maps and normal fields on $S.$ To be more precise, 
$\{(X_n)_{\sigma_n}\}_{n\in \n} \to X_\sigma$ in the  $\mathcal{C}^1(S)$-topology if $\{\|(X_n)_{\sigma_n}-X_\sigma\|_{1,S}\}_{n\in \n} \to 0.$
\end{definition}

Given $X_\sigma \subset \mathcal{M}^*(S),$ let $\partial_z X_\sigma=(\hat{\phi}_j)_{j=1,2,3}$ be the complex vectorial  1-form  of type $(1,0)$  on $S$ given by  $\partial_z X_\sigma|_M=\partial_z(X|_M),$ $\partial_z X_\sigma (\beta_{j}'(s)):= (X\circ \beta_{j})'(s) + i\sigma(\beta_{j}(s)) \wedge  (X\circ \beta_{j})'(s) ,$ 
where $s$ is the  arc length parameter of $X\circ \beta_{j}.$  
To be more precise, if $(U,z=x+i y)$ is a conformal chart on $\nb$ such that $\beta_j \cap U=z^{-1}(\r \cap z(U)),$ then $(\partial_z X_\sigma)|_{\beta_j \cap U}=\big[(X\circ \beta_j)'(s) + i   \sigma(\beta_{j}(s)) \wedge  (X\circ \beta_{j})'(s) \big]s'(x)dz|_{\beta \cap U},$ $j=1,\ldots,m.$ 

The analyticity of $\beta$ and the conformality property are crucial for the well-definition of $\partial_z X_\sigma$ on $\beta.$ In particular, $(\partial_z X_\sigma)|_{\beta_j}(P)=(\partial_z X_\sigma)|_{\beta_i}(P)$  at any point $P\in \beta_i\cap\beta_j,$  $i,j\in \{1,\ldots,m\}.$  As a consequence,  $\partial_z X_\sigma$  lies in 
 $\Wg_h(S)^3.$   
 
 Notice that $\sum_{j=1}^3 \hat{\phi_j}^2=0$ and   set $\hat{\eta}_1=\hat{\phi}_1-i \hat{\phi}_2,$ $\hat{\eta}_2=-\hat{\phi}_1-i \hat{\phi}_2$ and $\hat{g}=\hat{\eta}_2/\hat{\phi}_3.$ Since  $\hat{g}:S\to \overline{\c}$ is the stereographic projection of the generalized Gauss map $N_\sigma$ of $X_\sigma,$ it is conformal as well.

 \begin{remark}\label{re:conformal}
  $\hat{\phi}_j$ is a smooth object  on $S$ in the sense that  $\hat{\phi}_j/\theta$ is a smooth function, where $\theta$ is any never vanishing holomorphic 1-form on $\nb,$ $j=1,2,3.$ The same holds for $\hat{\eta}_i,$ $i=1,2.$  
  
In a similar way $\hat{g}\in \Fg_m(M).$ Furthermore, accordingly with Definition \ref{def:spacfun},  $\hat{g}$  lies in $\Fg_m(S)$  provided that $\hat{g} \neq \infty$ on $\beta-M.$
\end{remark}
 
Notice that $\sum_{j=1}^3 \hat{\phi}_j^2=0$ and $\mbox{Re} (\hat{\phi}_j)$ is an  exact  real 1-form on $S,$ $j=1,2,3.$ If $S$ is connected,  we also have $X(P)=X(Q)+\mbox{Re} \int_{Q}^P (\hat{\phi}_j)_{j=1,2,3},$ $P,$ $Q \in S.$  The pairs $(\hat{g},\hat{\phi}_3)$ and $(\hat{\eta}_j|_{M})_{j=1,2}$ will be called as the generalized  Weierstrass data  and  spinorial representation  of $X_\sigma,$ respectively. 
As $X|_M \in \mathcal{M}(M),$  then $(\phi_j)_{j=1,2,3}:=(\hat{\phi}_j|_{M})_{j=1,2,3},$ $(\eta_j)_{j=1,2}=(\hat{\eta}_j|_{M})_{j=1,2}$ and $g:=\hat{g}|_{M}$ are  the  Weierstrass data, spinorial representation and  meromorphic Gauss map of $X|_{M},$ respectively. Recall that all these data extend meromorphically to $M^c.$ 

The group homomorphism  $p_{X_\sigma}:\mathcal{H}_1(S,\z) \to \r^3, \quad p_{X_\sigma}(\gamma)=\mbox{Im} \int_\gamma \partial_z X_\sigma,$ is said to be the {\em generalized flux map} of $X_\sigma.$  Two marked immersions $X_{\sigma_1},$ $Y_{\sigma_2} \in \mathcal{M}^*(S)$ are said to be {\em flux equivalent} on $S$ if $p_{X_{\sigma_1}}=p_{Y_{\sigma_2}}.$

\begin{definition}\label{def:rest} 
Let $V\subset \nb$ be a finite collection of pairwise disjoint regions containing  $S,$  let $Y\in \mathcal{M}(V)$ and call $N:V\to\s^2$ its  Gauss map. We set 
$$\Rg_S(Y)=(Y|_{S},N|_{\beta}).$$ 
\end{definition}
Observe that $\Rg_S(\mathcal{M}(V))\subset\mathcal{M}^*(S),$ and notice that the {\em restriction map} $\Rg_S:\mathcal{M}(V) \to \mathcal{M}^*(S)$   is continuous with respect to the $\mathcal{C}^0(V)$-topology on $\mathcal{M}(V)$ and the $\mathcal{C}^1(S)$-topology on $\mathcal{M}^*(S).$ In the sequel, we write $$\text{$\|X _\sigma-Y\|_{1,S}:=\|X _\sigma-\Rg_S(Y))\|_{1,S}$ and $\|Z-Y\|_{1,S}:=\|\Rg_S(Z)-\Rg_S(Y))\|_{1,S}$}$$ for any $X_\sigma \in \mathcal{M}^*(S)$  and $Y,$ $Z\in \mathcal{M}(V).$ 

It is clear that   $p_{\Rg_S(Y)}=p_Y|_{\mathcal{H}_1(S,\z)}$ for any $Y \in \mathcal{M}(V),$ where $p_Y$ is the flux map of $Y.$

\begin{definition}
Given $X_\sigma \in \mathcal{M}^*(S),$ we denote by $\mathcal{M}_{X_\sigma}(\nb)$ the space of those immersions $Y \in \mathcal{M}(\nb)$ for which $\Rg_S(Y)$  is flux equivalent to $X_\sigma.$
\end{definition}

\section{Approximation by complete minimal surfaces with FTC} \label{sec:appro}

Roughly speaking, the aim of this section is to show that any finite collection of Jordan arcs and complete minimal surfaces with FTC and non-empty compact boundary (for instance, a finite collection of Jordan arcs and compact minimal surfaces), can  be uniformly approximated by connected complete minimal surfaces of FTC. Furthermore, the conformal structure and the flux map of the approximate sequence can be prescribed. This the message of the Fundamental Approximation Theorem below (see Theorem \ref{th:parabo} for a more general result).

Fix an open Riemann surface $\nb,$ and keep the notations of  Sections \ref{sec:pre} and \ref{sec:wei-spin}. Furthermore, assume that
\begin{itemize}
\item $\nb$ has {\em finite conformal type},
\item $S=M\cup \beta$ is an admissible subset in $\nb,$ and
\item $\nb-S$  consists of a finite collection of pairwise disjoint once punctured open discs.
\end{itemize}
In particular, $S$ is connected and $j_*:\mathcal{H}_1(S,\z) \to \mathcal{H}_1(\nb,\z)$ is an isomorphism, where $j:S\to \nb$ is the inclusion map.

Label by $\nu$ the genus of the Osserman compatification $\nb^c$ of $\nb,$ notice that $M^c\subset S^c\subset \nb^c.$ Put $M=M^c-\{E_{1},\ldots,E_{a}\},$  $\nb=\nb^c-\{E_{1},\ldots,E_{a+b}\}$ and $\nb_0=\nb^c-\{E_{a+1},\ldots,E_{a+b}\},$ for suitable points $E_1,\ldots,E_{a+b}\in \nb.$
Label  $U_1,\ldots,U_b$  as the connected components (open discs) of $N^c-S^c,$ where up to relabeling  $E_{a+j} \in {U_j},$ $j=1,\ldots,b.$

\begin{theorem}[The Fundamental Approximation Theorem] \label{th:density}

For any   $X_\sigma \in \mathcal{M}^*(S),$ there exists a sequence  $\{Y_n\}_{n \in \n}  \subset \mathcal{M}_{X_\sigma}(\nb)$  such that $\{\Rg_S(Y_n)\}_{n \in \n} \to X_\sigma$ in the ${\cal C}^1(S)$-topology. 

Furthermore, if $C$ is a positive constant and  $V$ a closed tubular  neighborhood of $S$ in $\nb,$    $\{Y_n\}_{n \in \n}$ can be chosen in such a way that $d_{Y_n}(S, \partial(V)) \geq C$ for all $n,$ where $d_{Y_n}$ is the intrinsic distance in $N$ induced by $Y_n.$
\end{theorem}

The global strategy for proving this theorem has essentially three phases.
\begin{enumerate}[(I)]
\item {\em First phase:} Show that the spinorial representation of $X_\sigma$ on $S$ can be approximated by  holomorphic spinorial data on $\nb$ extending meromorphically to $\nb^c$ (this technical result corresponds to Lemma \ref{lem:sp} in paragraph \ref{subsub:phase1}). 
\item {\em Second phase:} Prove that the  approximating sequence of meromorphic spinorial data on $\nb^c$ can be slightly deformed in order to solve the period problem (this part corresponds to Lemmas \ref{lem:implicit} and \ref{lem:4} in paragraph \ref{subsub:phase2}).
\item {\em Third phase:} Conclude the proof of Theorem  \ref{th:density} (see paragraph \ref{subsub:phase3}).
\end{enumerate}

\subsection{Proof of the Fundamental Approximation Theorem} \label{subsec:algbrid}
Before starting with the first phase of the program, we  establish some basic conventions that can be assumed without loss of generality. 
This is the content of the following three propositions.

Take $X_\sigma\in \mathcal{M}^*(S)$ as in the statement of Theorem \ref{th:density}. 

The first proposition  simply says  that $X(M)$ can be supposed without containing planar domains.  
\begin{proposition}\label{pro:nonflat} 
Without loss of generality, we can suppose that $X(M)$ contains no planar domains.
\end{proposition}
\begin{proof} Let us show that there exists  a sequence $\{Y^j_{\sigma_j}\}_{j \in \n}\subset \mathcal{M}^*(S)$ such that $\{Y^j_{\sigma_j}\}_{j \in \n} \to X_\sigma$ uniformly on $S,$ $Y^j_{\sigma_j}$ is flux equivalent to $X_\sigma$ on $S$ and $Y^j(M)$ contains no planar domains, $j \in \n.$

Indeed, since any flat minimal surface can be approximated by non flat ones, we can find $\{Y^j\}_{j \in \n}\subset {\cal M}(M)$ such that $Y^j(M)$ contains no planar domains for all $j$ and $\{Y^j\}_{j \in \n}\to X|_M$ in the ${\cal C}^0(M)$-topology. Write $N^j$ for the  Gauss map of $Y^j,$ and  extend $Y^j$ and $N_j$ to $\beta$ in a smooth and conformal way so that $(Y^j,N^j|_\beta)\in \mathcal{M}^*(S)$ and
$(Y^j,N^j|_\beta)$ is flux equivalent to $X_\sigma$ for all $j,$ and $(Y^j,N^j|_\beta)\}_{j \in \n} \to X_\sigma$ in the ${\cal C}^1(S)$-topology. 

To finish, notice that if the Fundamental Approximation Theorem held in the non-flat case, the immersions $(Y^j,N^j|_\beta)$ would lie in the closure of $\Rg_S(\mathcal{M}_{X_\sigma}(\nb))$ in $\mathcal{M}^*(S),$ were $\Rg_S$ is the restriction map in Definition \ref{def:rest}, hence the same would occur for $X$ and the first part of the theorem would hold. The second one can also be guaranteed in the process. 
\end{proof}

Label $\partial_z X_\sigma=(\hat{\phi}_j)_{j=1,2,3},$ and consider the generalized  Weierstrass data  $(\hat{g},\hat{\phi}_3)$ and  spinorial representation  $(\hat{\eta}_j|_{M})_{j=1,2}$ of $X_\sigma.$  Write $d \hat{g}$ for the  1-form  of type $(1,0)$ on $S^c$  given by  $d \hat{g}|_{M^c}=d (\hat{g}|_{M^c})$ and  
$d\hat{g}(\beta'(s))=(\hat{g}\circ \beta)'(s),$ where $s$ is the arc length parameter of $X\circ \beta.$ In other words, if $(U,z=x+i y)$ is a conformal chart in $N$ so that $\beta \cap U=z^{-1}(\r \cap z(U)),$ then $d\hat{g}|_{\beta \cap U}=(\hat{g}\circ \beta)'(s) s'(x)dz|_{\beta \cap U}.$ Since $\hat{g}$ is conformal (see Remark \ref{re:conformal}), it is not hard to check that $d\hat{g}$ is well defined. Furthermore,  $d\hat{g} \in \Wg_m(S)$ when $\hat{g}(P) \neq \infty$ for any $P \in \beta.$  
 Write $(\phi_j)_{j=1,2,3}=(\hat{\phi}_j|_{M})_{j=1,2,3},$ $(\eta_j)_{j=1,2}=(\hat{\eta}_j|_{M})_{j=1,2}$ and $g=\hat{g}|_{M}$ for  the  Weierstrass data, the spinorial representation and the meromorphic Gauss map of $X|_{M},$ respectively, and call with the same name their meromorphic extensions to $M^c.$

The second convention deals with the behavior of $\hat{g}$ on $\partial(S)=\partial(M) \cup \beta.$

\begin{proposition} \label{pro:motion}
Without loss of generality,  we can assume that  
\begin{enumerate}[(i)]
\item $\hat{g},$ $1/\hat{g},$ $(\hat{g}^2-1),$ and $d \hat{g}$  never vanish on $\partial(S),$ hence the same holds for $\hat{\eta}_i,$  $i=1,2,$ and $\hat{\phi}_j,$ $j=1,2,3$ (in particular, $\hat{g}  \in \Fg_m(S)$ and  $d \hat{g} \in \Wg_m(S)$),
\item $d\hat{g}\neq 0$ at any point of $\hat{g}^{-1}(\{0,\infty\}),$ and
\item $\hat{g}(E_i)\neq 0,\infty,$ $i=1,\ldots,a.$ 
\end{enumerate}
In particular, $m_i:=\mbox{Ord}_{E_i}(\hat{\phi}_3)=\mbox{Ord}_{E_i}(\hat{\eta}_1)=\mbox{Ord}_{E_i}(\hat{\eta}_2)>1,$ where $\mbox{Ord}_{E_i}(\cdot)$ means pole order at $E_i,$ $i=1,\ldots,a.$

\end{proposition}
\begin{proof} Up to a rigid motion, we can suppose that $\hat{g}(E_i)\neq 0,\infty,$ $i=1,\ldots,a,$ and  $d\hat{g}\neq 0$ at any point of $\hat{g}^{-1}(\{0,\infty\})\cap M.$ In particular, $\mbox{Ord}_{E_i}(\hat{\phi}_3)=\mbox{Ord}_{E_i}(\hat{\eta}_1)=\mbox{Ord}_{E_i}(\hat{\eta}_2)>1,$ $i=1,\ldots,a.$

Recall that $X|_M$ is non flat and extends  as a conformal minimal immersion beyond $M$ in $\nb.$ Therefore, we can find a sequence $M_{(1)} \supset M_{(2)} \supset \ldots $ of closed tubular  neighborhoods of $M$ in $\nb$ such that $M_{(j)} \subset M_{(j-1)}^\circ$ for any $j,$ $M=\cap_{j \in  \n} M_{(j)},$ $X$ and $\hat{g}$ extend  (with the same name) as a conformal minimal immersion and a meromorphic function to $M_{(j)},$  $\hat{g},$ $1/\hat{g},$ $(\hat{g}^2-1),$ and $d \hat{g}$  never vanish on $\partial({M}_{(j)})$ for all $j,$  and $d\hat{g}\neq 0$ at any point of $\hat{g}^{-1}(\{0,\infty\})\cap M_{(j)}$ for all $j.$ Call  $\beta_{(j)}:=\beta-M_{(j)}^\circ,$ and without loss of generality assume that $S_j:=M_{(j)}\cup \beta_{(j)}$ is admissible in $\nb$ as well for all $j.$

Up to suitably deforming $X|_\beta$ and $\sigma|_\beta,$ we can construct marked immersions  $Z^j_{\sigma_j}\in \mathcal{M}^*(S_j),$ $j\in \n,$ such that
\begin{itemize}
\item  $Z^j|_{M_{(j)}} =X|_{M_{(j)}}$ and $Z^j|_{M_{(j)}}$ is flux equivalent to $X_\sigma$ on $S,$
\item $\hat{g}_j,$ $1/\hat{g}_j,$ $(\hat{g}_j^2-1),$ and $d \hat{g}_j \neq 0$ on $\partial(S_j),$ and $d\hat{g}_j\neq 0$ at any point of $\hat{g}_j^{-1}(\{0,\infty\})\cap M_{(j)},$ where $\hat{g}_j$ is the generalized  Gauss map  of $Z^j,$ and 
\item   $\{(Z^j|_S,N_{\sigma_j}|_\beta)\}_{j \in \n} \to X_\sigma$ in the ${\cal C}^1(S)$-topology, where $N_{\sigma_j}$ is the Gauss map of $Z^j_{\sigma_j}.$
\end{itemize}

If Theorem \ref{th:density} held for $Z^j_{\sigma_j},$ $j \in \n,$ we would infer that  $Z^j_{\sigma_j}$ lies in the closure of $\Rg_{S_j}(\mathcal{M}_{X_\sigma}(\nb))$ in $\mathcal{M}^*(S_j),$ $j \in \n.$ Since $\{(Z^j|_S,N_{\sigma_j}|_\beta)\}_{j \in \n} \to X_\sigma$ in the ${\cal C}^1(S)$-topology, we would infer that  $X_\sigma$ lies in the closure of $\Rg_{S}(\mathcal{M}_{X_\sigma}(\nb))$ in $\mathcal{M}^*(S)$ as well and we are done.

 The second part of the theorem would also be achieved in the process. 
 
\end{proof}

Let us go to the first phase of the program.

\subsubsection{Approximating the spinorial data of $X_\sigma$ on $S$ by global holomorphic ones in $\nb.$} \label{subsub:phase1}

The following notation is previously required.

Let $\Theta_j$ denote the class of $\eta_j$ in $\frac{\Sg_m(M)}{\sim},$ and for the sake of simplicity, write  $\xi_j$ for the associated morphism $\xi_{\Theta_j}:\mathcal{H}_1(M,\z) \to \z_2$ (notice that these objects make sense even when $M$ is  not connected). Let us show that there is a canonical extension of $\xi_j$ to $\mathcal{H}_1(\nb,\z)$ depending on $\hat{\eta}_j.$ Indeed, recall that $\mathcal{H}_1(\nb,\z)=\mathcal{H}_1(S,\z)$ and take an arbitrary closed curve $c\in \mathcal{H}_1(S,\z)$. From Proposition \ref{pro:motion},  all the zeros of $\eta_j=\hat{\eta}_j|_M$ have even order and  $\hat{\eta}_j$ never vanishes on $\partial(S).$ If we take any conformal annulus $(A,z)$ in $\nb$ such that $A$ is a closed tubular  neighborhood  of $c,$ it suffices to set  $\xi_j(c)=0$ when $\sqrt{\hat{\eta}_j(z)/dz}$  has  a well defined branch along $c$ and  $\xi_s(c)=1$ otherwise (this computation does not depend on the chosen $(A,z)$).  

On the other hand, the fact that $\hat{\eta}_2/\hat{\eta}_1=\hat{g}^2$ implies that  $\xi_1=\xi_2,$ hence one can say that $\hat{\eta}_1$ and $\hat{\eta}_2$   are "spinorially equivalent" on $S.$  Lemma \ref{lem:spin} guarantees the existence of a unique spinor structure on $\nb$ associated to $\xi_1.$  By definition, an 1-form $\theta \in \Sg_m(\nb)$ is said to be spinorially equivalent to $\hat{\eta}_1$ (and so to $\hat{\eta}_2$)  if $\xi_\Theta=\xi_1,$ where $\Theta \in \frac{\Sg_m(\nb)}{\sim}$ is the class of $\theta.$   By  Lemma \ref{lem:spin}, we can always find 1-forms of this kind in  $\Sg_m(\nb).$\\

The main goal of this phase is to prove the following:

\begin{lemma} \label{lem:sp}
 There are $\{\eta_1^n\}_{n \in \n},$ $\{\eta_2^n\}_{n \in \n} \subset \Sg_h(\nb)$ such that:
\begin{enumerate}[(i)]
  \item $\{\eta_j^n|_{S}\}_{n \in \n} \to \hat{\eta}_j$ in the $\mathcal{C}^0(S)$-topology,  
$\eta_j^n$ never vanishes on $\partial(S),$  $(\eta_j^n|_{M^c})=(\hat{\eta}_j|_{M^c}),$    $(\eta_j^n|_{M^c}-\hat{\eta}_j|_{M^c})_0 \geq \prod_{i=1}^a E_i^{m_i},$ and  
$(\eta_j^n)_\infty \geq \prod_{k=a+1}^{a+b} E_k,$ $j\in \{1,2\},$  $n \in \n.$ 
  \item $\eta_1^n$ and $\eta_2^n$ are spinorially equivalent in $\Sg_m(\nb)$ and have no common zeros on $\nb.$
 \end{enumerate}   
\end{lemma}
\begin{proof} The following claims will be useful:

\begin{quote} 
\begin{assertion} \label{ass:Jacobi}
Let $R$ be a compact Riemann surface with empty boundary and genus $\nu.$ Given an open disc $U \subset R,$ a point $Q \in R$ and a divisor $D_1 \in \div(R),$ there exists an integral divisor $D_2 \in \div(U)$ of degree $\nu$ and $n_0 \in \n$  such that $D_2^{n_0}D_1^{-1} Q^{-n_1}$ is the principal divisor associated to some $f \in \Fg_m(R),$ where $n_1={n_0} \nu-\deg(D_1).$
\end{assertion}
\end{quote}
\begin{proof} Since the proof is trivial when $\nu=0,$ we will assume that $\nu \geq 1.$ By Jacobi's theorem, we can find an open disc $W \subset U$  such that $\varphi_Q:W_\nu \to \varphi_Q(W_\nu)$ is a diffeomorphism, where $\varphi_Q$ is the Abel-Jacobi map with base point $Q$ and $W_\nu$ is the set of divisors in $R_\nu$  with support in $U.$ Since $J(R)$ is a compact additive Lie Group and $\varphi_Q(W_\nu) \subset J(R)$ is an open subset, for large enough ${n_0} \in \n$ one has ${n_0} \varphi_E(W_\nu)=J(R).$ Therefore, there is $D_2 \in W_\nu$ such that $\varphi_Q(D_2^{n_0})=\varphi_Q(D_1)=\varphi_Q(D_1 Q^{n_1}),$ where $n_1={n_0} \nu -\deg(D_1).$ The claim follows from Abel's theorem. 
\end{proof}

\begin{quote}
\begin{assertion} \label{ass:otra} We can find $\theta_1,$ $\theta_2 \in \Sg_h(\nb)$ so that $|\theta_1|+|\theta_2|$ has no zeros in $\nb,$  $\theta_j$ is spinorially equivalent to $\hat{\eta}_1$ and $\hat{\eta}_2,$  $\theta_j$  never vanishes on $S$ and $(\theta_j)_\infty \geq \prod_{i=1}^a E_i^{2 m_i},$  $j=1,2.$
\end{assertion}
\end{quote}
\begin{proof} Take $\theta \in \Sg_m(\nb)$ spinorially equivalent to $\hat{\eta}_1$ and  $\hat{\eta}_2$ (see Lemma \ref{lem:spin}).

Let $k_i$ denote the zero order of $\theta$ at $E_i$ ($k_i=0$ provided that $\theta(E_i) \neq 0$), $i=1,\ldots,a,$ write $(\theta|_{\nb})=D_0^2,$ and fix two disjoint open discs $V_1,$ $V_2 \subset \nb^c-(S^c).$ By Claim \ref{ass:Jacobi},  there are $D_j \in \div(V_j)$ of degree $\nu,$ $n_j \in \n$ and  $h_j \in \Fg_m(\nb^c)$ such that $(h_j)=D_j^{n_j} E_{a+b}^{-v_j} D_0^{-1} \prod_{i=1}^a E_i^{-m_i-k_i},$ where  $v_j=n_j \nu-\deg(D_0)-\sum_{i=1}^a (m_i+k_i),$ $j=1,2.$ It suffices to put $\theta_j=h_j^2 \theta,$ $j=1,2.$
\end{proof}

Let $m_{i,j}\geq 2 m_i$  denote the pole order of $\theta_j$ at $E_i,$ $i=1,\ldots,a,$ and likewise call $n_{k,j}$ as the zero order of $\theta_j$ at $E_{a+k}$ ($n_{k,j}=0$ provided that $\theta_j(E_{a+k})\neq 0$),   $k=1,\ldots,b.$ Set $s_j=\frac{\hat{\eta}_j}{\theta_j},$ $j=1,2,$  and observe that $ s_j \in \Fg_h(S^c).$ Moreover, Proposition \ref{pro:motion} and Claim \ref{ass:otra} give that  $s_j\neq 0,$ $\infty$ on $\partial({S}),$ $(s_j|_{M^c})=(\hat{\eta}_j)_0 \prod_{i=1}^a E_i^{m_{i,j}-m_i} \geq \prod_{i=1}^a E_i^{m_i},$ $j=1,2,$ and $|s_1|+|s_2|$ has no zeros in $S.$ Claim \ref{ass:otra} also says that   $s_j=t_j^2$ for some $t_j \in \Fg_h(S^c),$ $j=1,2.$\\

Let us construct $\eta_1^n.$  

Consider  a collection $C_1$  of pairwise disjoint closed discs  in $\nb_0-S^c$ containing all  the zeros of $\theta_2|_{\nb}$ (recall that $|\theta_2|$ never vanishes on $S^c,$ see Claim \ref{ass:otra}) and meeting  all the bounded components of $\nb^c-S^c.$ It is clear that $S^c \cup C_1$ is admissible in the open Riemann surface $\nb_0=S^c\cup \nb.$ Consider the continuous map  ${t}_1^{*}:S^c \cup C_1\cup D_1 \to \c,$  ${t}_1^{*}|_{S^c}=t_1,$ ${t}_1^{*}|_{C_1}=\delta,$ where $\delta$ is a non-zero constant, and notice that $t_1^* \in \Fg_h(S^c \cup C_1).$  By Theorem \ref{th:runge} applied to the open Riemann surface $\nb_0,$ the admissible subset $S^c\cup C_1,$ the  function ${t}_1^{*}\in \Fg_h(S^c \cup C_1),$ and the divisor $(\hat{\eta}_j|_{M^c})_0 \prod_{i=1}^a E_i^{m_{i,1}+m_i},$ we can find  $\{H_{n,1}\}_{n\in \n} \subset \Fg_h(\nb_0)\subset \Fg_m(\nb^c)$  such that $\{H_{n,1}-{t}_1^{*}\}_{n\in \n} \to 0$ uniformly on $S^c \cup C_1$ and  $(H_{n,1}|_{M^c}-t_1|_{M^c})_0 \geq (\hat{\eta}_j|_{M^c})_0 \prod_{i=1}^a E_i^{m_{i,1}+m_i}.$ In particular, $(H_{n,1}/t_1)|_{M^c}$ is holomorphic and $(H_{n,1}/t_1)|_{S^c}$ and $H_{n,1}|_{C_1}$ are  never-vanishing for large enough $n$ (without loss of generality,  for all $n$), see Remark \ref{re:runge}.

\begin{assertion} \label{ass:infinito}
Without loss of generality, we can  assume that the sequence of pole multiplicities $$\{\mbox{Ord}_{E_{a+k}} (H_{n,1})\}_{n\in \n}$$ is divergent for all $k=1,\ldots,b.$ In particular,  we can assume that   $(H_{n,1})_\infty \geq \prod_{k=1}^b E_{a+k}^{(n_{k,1}+1)/2}$ for all $n\in \n.$  
\end{assertion}
\begin{proof} From Riemann-Roch theorem, it is not hard to find a function $T\in \Fg_h(\nb_0)\cap \Fg_m(\nb)$ such that  $(T)_\infty\geq \prod_{k=1}^b E_{a+k}$ and $(T)_0\geq (\hat{\eta}_j|_{M^c})_0 \prod_{i=1}^a E_i^{m_{i,1}+m_i}.$ 

For each $n\in \n,$ take $j_n\in \n$ such that $(H_{n,1} T^{j_n})_\infty\geq \prod_{k=1}^b E_{a+k}^{n_{k,1}+n},$ and then choose   $k_n\in \n$ such that $|H_{n,1}T^{j_n}|,$ $|T^{j_n}|<k_n/n$ on $S^c\cup C_1.$ 

The sequence  $\{H_{n,1} (\frac{T^{j_n}}{k_n}+1)\}_{n\in \n}$ formally satisfies the same properties as $\{H_{n,1}\}_{n\in \n}$ and has the desired pole orders. To finish, replace $H_{n,1}$ for $H_{n,1} (\frac{T^{j_n}}{k_n}+1)$ for all $n.$

\end{proof}

Call  $F_{n,1}= ({H_{n,1}})^2,$ and notice that $\{F_{n,1}\}_{n \in \n} \to s_1$ in the 
$\mathcal{C}^0(S)$-topology, $F_{n,1}/s_1$  never vanishes on $S^c$ and is holomorphic on $M^c,$  $F_{n,1}$ never vanishes on $\partial(S),$ $(F_{n,1}\theta_1 |_{M^c})=(\hat{\eta}_1|_{M^c}),$  $\big((F_{n,1}-s_1)|_{M^c}\big)_0 \geq (\hat{\eta}_1|_{M^c})_0\prod_{i=1}^a E_i^{m_{i,1}+m_i}$ and $(F_{n,1})_\infty \geq \prod_{k=a+1}^{a+b} E_k^{n_{k,1}+1}$ for all $n.$ 

Setting $\eta_1^n:= F_{n,1} \theta_1,$ item $(i)$ holds for $j=1.$ \\

For constructing $\eta_2^n,$ we reason in a similar way. 

Choose a collection $C_{n,2}$ of pairwise disjoint closed discs  in $\nb_0-S^c$  containing all the zeros of $\eta_1^n$ in $\nb_0-M^c$ and meeting all the bounded components of $\nb^c-S^c.$ Set ${t}_{n,2}^{*}:S^c \cup C_{n,2} \to \c,$  ${t}_{n,2}^{*}|_{S^c}=t_2$ and ${t}_{n,2}^{*}|_{C_{n,2}}=\delta,$ 
  where $\delta$ is any non-zero constant.  As above we can construct $H_{n,2} \in \Fg_m(\nb^c) \cap {\Fg_h}(N_0)$ satisfying that $|H_{n,2}-{t}_{n,2}^{*}|<1/n$ on $S^c \cup C_{n,2},$  $(H_{n,2}|_{M^c}-t_2|_{M^c})_0 \geq (\hat{\eta}_2|{M^c})_0 \prod_{i=1}^a E_i^{m_{i,2}+m_i},$  $H_{n,2}/t_2$ is holomorphic and never-vanishing on $S^c,$ $H_{n,2}|_{C_{n,2}}$ is holomorphic and never-vanishing,  and $\mbox{Ord}_{E_k} H_{n,2}\geq (n_{k,2}+1)/2$ for all $n\in \n,$ $k=1,\ldots,b.$ Set  $F_{n,2}= (H_{n,2})^2,$ and observe that $\{F_{n,2}\}_{n \in \n} \to s_2$ in the $\mathcal{C}^0(S)$-topology,  the holomorphic function $F_{n,2}/s_2$ never vanishes on $S^c$ for all $n,$  $(F_{n,2}\theta_2 |_{M^c})=(\hat{\eta}_2|_{M^c}),$  $\big((F_{n,2}-s_2)|_{M^c}\big)_0 \geq (\hat{\eta}_2|_{M^c})_0\prod_{i=1}^a E_i^{m_{i,2}+m_i}$ and $(F_{n,2})_\infty \geq \prod_{k=a+1}^{a+b} E_k^{n_{k,2}+1}$ for all $n.$

Choosing $\eta_2^n:= F_{n,2} \theta_2,$  item $(i)$ holds for $j=2.$\\ 

Finally, let's check item $(ii).$ Obviously, $\eta_1^n$ and $\eta_2^n$ are spinorially equivalent in $\Sg_m(N).$ Recall that $F_{n,1}$, $\theta_1 \cdot \theta_2,$ and $|\theta_1|+|\theta_2|$ never vanish on $C_1,$ $S$ and $\nb,$ respectively. Therefore, from the choice of $C_1,$ one has that  $|\eta_1^n|+|\theta_2|$ has no zeros on $\nb,$ $n \in \n.$ Likewise, the choice of $C_{n,2}$ and the fact that $F_{n,2}|_{C_{n,2}}$  never vanishes imply that $|\eta_1^n|+| \eta_2^n|$ have no zeros on $N-S.$ Moreover, $\eta_1^n$ and $\eta_2^n$ never vanishes on $\partial(S)$ and $(\eta_j^n |_{M^c})=(\hat{\eta}_j),$ $j=1,2,$ hence  $|\eta_1^n|+|\eta_2^n|$ never vanishes on  ${M}$ as well and we are done. 
\end{proof}
\begin{remark} \label{re:data}
By (ii) in Lemma \ref{lem:sp}, the 1-form $\phi_3^n:= \sqrt{\eta_1^n \eta_2^n}$ is well defined and lies in $\Wg_m(\nb^c) \cap \Wg_h(\nb).$ With the proper choice of the square root branch,  $\phi_3^n|_{{M}}-\phi_3$ extends holomorphically to $M^c$ and $\{\phi_3^n|_{S^c}-\hat{\phi}_3\}_{n \in \n} \to 0$ in the $\mathcal{C}^0(S)$-topology. 
In other words, if we call  $\Phi_n:=(\phi_j^n)_{j=1,2,3} \in \Wg_h(\nb)^3$ the Weierstrass data associated to $(\eta_1^n,\eta_2^n)$ by equation \eqref{eq:spinor}, $n\in \n,$  then $\{\Phi_n|_S\}_{n\in \n}$  converge  in the $\mathcal{C}^0(S)$-topology to $(\hat{\phi}_j)_{j=1,2,3}.$
\end{remark}


The first tentative of solution for the Fundamental Approximation Theorem could be to choose $Y_n:=\mbox{Re}\int \Phi_n,$ $n\in \n.$ However,  $\{\Phi_n\}_{n\in \n}$ may have real periods (the immersions $\{Y_n\}_{n\in \n}$ could not be well defined), and we have no control on the associated flux maps. 

At this point we  start with the second phase of the program.

\subsubsection{Deforming the global spinorial data and solving the period problem.} \label{subsub:phase2}

In order to overcome the above problems, it is  necessary to slightly deform these data in a suitable way. 

We need the following
\begin{definition}
Fix a homology basis $B_0$ of  $\mathcal{H}_1(S^c,\z),$ hence of $\mathcal{H}_1(\nb_0,\z),$ and call $\varsigma_0=3(2\nu+b-1)$ the cardinal number of $B_0.$
\end{definition}

Roughly speaking, our global strategy consists of the following.\\

Firstly, we present the natural space of deformations. In our case, it corresponds to 
$${\cal L}=\{f \in \Fg_h(\nb_0) \;:\; (f)_0 \geq \prod_{j=1}^a E_j^{m_i}\}\subset \Fg_m(\nb^c).$$ 
By Riemann-Roch theorem, ${\cal L}$ is a linear subspace of $\Fg_m(\nb^c)$ with  infinite dimension and finite codimension. Up to restriction to $S,$ ${\cal L}$ can be viewed as subspace of the  complex normed space $(\Fg_h(S^c),\| \cdot \,\|_{0,S}),$ where  $\|h\|_{0,S}=\max_{S^c} |h|=\max_{S} |h|$ is the norm of the maximum on $S^c.$

Then, we introduce an analytical deformation $\{\hat{\Phi}(f),|\;f\in {\cal L}\}\subset \Wg_h(S)^3$ of  $\hat{\Phi},$   where  $\hat{\Phi}({\bf 0})=\hat{\Phi}$ (here ${\bf 0}$ is the constant zero function), and likewise for $\Phi_n,$ $n\in \n.$  Subsequently, we define the Fréchet differentiable analytical period operators
$${\cal P}:{\cal L} \to \c^{\varsigma_0}, \quad {\cal P}(f)= \big(\int_{d} \hat{\Phi}(f) -\hat{\Phi}\big)_{d  \in B_0},$$
$${\cal P}_n:{\cal L} \to \c^{\varsigma_0}, \quad {\cal P}_n(f)= \big(\int_{d} \Phi_n(f) -\hat{\Phi}\big)_{d  \in B_0},\quad n \in \n.$$
The key step is to prove that $d {\cal P}_0$ is surjective (Lemma \ref{lem:implicit} below), and consequently that ${\cal P}({\bf 0})=0\in \c^{\varsigma_0}$ is an interior point of ${\cal P}({\cal L}).$ Since $\{{\cal P}_n\}_{n\in \n}\to {\cal P}$ in a uniform way,  we can deduce that $0$ is an interior point of ${\cal P}_n({\cal L})$ as well, $n$ large enough, and  find  ${\bf h}_n \in {\cal P}_n^{-1}(0),$ $n\in \n,$ such that $\{{\bf h}_n\}_{n\in\n}\to {\bf 0}.$  

Therefore, the sequence $\{\Phi_n({\bf h}_n)\}_{n\in \n}$ uniformly approximates $\hat{\Phi}$ on $S,$ have no real periods on $S,$ and induce the same flux map as $\hat{\Phi},$ concluding the second phase.\\

Let us develop carefully this program.

For  each $f \in {\cal L},$ set $\hat{\eta}_j(f)=(1+j f)^2 \hat{\eta}_j,$ $j=1,2,$ define $\hat{\phi}_k(f)$ following equation \eqref{eq:spinor},  and notice that  $\hat{\phi}_k(f)-\hat{\phi}_k\in \Fg_h(S^c),$ $k=1,2,3.$  

Endow ${\cal L}$ with the norm $\|\cdot\|_{0,S}$ of the maximum on $S$  inducing the $\mathcal{C}^0(S)$-topology. By the maximum principle, this norm  coincides with the one of the maximum on  $S^c.$ 
Consider  the Fréchet differentiable map
$${\cal P}:{\cal L} \to \c^{\varsigma_0}, \;{\cal P}(f)=[ \big(\int_{d} \hat{\phi}_j(f) -\hat{\phi}_j\big)_{d  \in B_0}]_{j=1,2,3}.$$  It is clear that ${\cal P}({\bf 0})=0\in \c^{\varsigma_0},$ where ${\bf 0}$ is the constant zero function.

\begin{lemma} \label{lem:implicit}
The complex Fréchet derivative $d {\cal P}|_0:{\cal L} \to \c^{\varsigma_0}$ of ${\cal P}$ at ${\bf 0}$ is surjective.
\end{lemma}
\begin{proof} 

Reason by contradiction, and assume that  $d {\cal P}|_0({\cal L})$ lies in a hyperplane $U= \{((x_d^j)_{d  \in B_0})_{j=1,2,3} \in \c^{\varsigma_0}\;:\; \sum_{j=1}^3\big(\sum_{d  \in B_0} \lambda_d^j x_d^j \big)=0\},$ where $\sum_{j=1}^3\big(\sum_{d  \in B_0}|\lambda_d^j| \big)\neq 0.$ 

Therefore $d {\cal P}|_0(f)=\frac{d{\cal P}(t f)}{dt}|_{_{t=0}}  \in U,$ for any $f \in {\cal L},$ that is to say
\begin{equation} \label{eq:fun}
\int_{\Gamma_1} f \hat{\eta}_1 + \int_{\Gamma_2} f \hat{\eta}_2+ \int_{\Gamma_3} f \hat{\phi}_3=0, \;  \mbox{for all}\; f \in {\cal L},
\end{equation}
where $\Gamma_j \in \mathcal{H}_1(S^c,\c),$ $j=1,2,3,$ are the cycles with complex coefficients given by:  $$\Gamma_1=\sum_{d  \in B_0} (\lambda_d^1+i \lambda_d^2)\, d, \; \Gamma_2=2\sum_{d  \in B_0} (-\lambda_d^1+i \lambda_d^2) \,d,\; \Gamma_3=3\sum_{d  \in B_0} \lambda_d^3 \,d.$$ 

The idea of the proof is to show that equation \eqref{eq:fun} yields that $\Gamma_1=\Gamma_2=\Gamma_3=0,$ a contradiction. Let us go to the details.

From Proposition \ref{pro:motion}-(i),  $\hat{\phi}_k,$ $\hat{\eta}_j\in \Wg_h(S),$ $d \hat{g}\in \Wg_m(S),$ and $\hat{g}\in\Fg_m(S)$  are  never vanishing objects on $\partial(S),$ and therefore, their associated divisors have support in $M^c-\partial(M).$  This fact is crucial for a good understanding of the following notations and arguments.

Set ${\cal L}_0=\{f \in \Fg_h(\nb_0) \;:\; (f)_0 \geq (\hat{\phi}_3)_0^2\prod_{i=1}^a E_i^{m_i}=(\hat{g})_0^2 (\hat{g})_\infty^2\prod_{i=1}^a E_i^{m_i}\}\subset {\cal L}.$ From Riemann-Roch theorem, ${\cal L}_0$ is a linear subspace of ${\cal L}$ of infinite dimension. 
Since $m_i$ is the pole order of $\phi_3$ at $E_i$ and $g(E_i) \neq 0,$ $\infty$ for all $i=1,\ldots,a$ (see Proposition \ref{pro:motion}-(iii)), then $df/\hat{\phi}_3\in \Fg_h(S^c)$ and $\left( (d f|_{M^c})/\hat{\phi}_3\right) \geq \prod_{i=1}^a E_i^{m_i}$ for all $f \in {\cal L}_0.$ By Theorem \ref{th:runge} and Remark \ref{re:runge}, there is $\{f_n\}_{n \in \n} \subset {\cal L}$ converging to $df/\hat{\phi}_3$ in the $\mathcal{C}^0(S)$-topology. Applying  equation (\ref{eq:fun}) to $f_n$ and taking the limit as $n$ goes to $+\infty,$ we infer that  $ \int_{\Gamma_1} \frac{df}{\hat{g}}  + \int_{\Gamma_2} \hat{g} df=0$ for any $f \in {\cal L}_0.$  Integrating by parts,
\begin{equation} \label{eq:fun1}
\int_{\Gamma_1} \frac{f d \hat{g}}{\hat{g}^2}-\int_{\Gamma_2} f d \hat{g}=0,\quad \mbox{for all} \;f \in {\cal L}_0.
\end{equation}

Denote by  ${\cal L}_{1}=\{f \in \Fg_h(\nb_0) \,:\; (f)_0 \geq ({\hat{g}}^2-1)_0^2 (d {\hat{g}})_0^2 \prod_{i=1}^a E_i^{2m_i} \}\subset {\cal L}.$  As above, from Riemann-Roch theorem  ${\cal L}_{1}$ is a linear subspace of ${\cal L}$ of  infinite dimension. For any $f \in {\cal L}_{1},$ the function  $h_{f}:=\frac{\hat{g}^2 df}{(\hat{g}^2-1) d \hat{g}}$ lies in $\Fg_h(S^c)$ and satisfies that $(h_{f}) \geq (\hat{\phi}_3)_0^2  \prod_{i=1}^a E_i^{m_i}$ (take into account Proposition \ref{pro:motion}-(ii)). By Theorem \ref{th:runge} and Remark \ref{re:runge}, $h_{f}$ lies in the closure of ${\cal L}_0$ in $(\Fg_h(S^c),\|\cdot\|_{0,S}),$ hence equation (\ref{eq:fun1}) can be formally applied  to $h_{f}$ to obtain that $\int_{\Gamma_1-\Gamma_2} \frac{df}{\hat{g}^2-1}=0,$ for any $f \in {\cal L}_{1}.$ Integrating by parts,
\begin{equation} \label{eq:fun2}
 \int_{\Gamma_1-\Gamma_2}  \frac{f d \hat{g}}{(\hat{g}^2-1)^2}=0,\quad \mbox{for all} \; f \in {\cal L}_{1}.
\end{equation}

At this point, we need the following 
\begin{assertion} \label{ass:puntos}
For any  $P_1,\ldots,P_r\in \nb_0,$  $n_1,\ldots,n_r\in \n,$ and  $\tau \in \Wg_h(\nb_0),$  there exists $F\in \Fg_h(\nb_0)$ such hat $(\tau+dF)_0 \geq \prod_{j=1}^r P_j^{n_j}.$
\end{assertion}
\begin{proof} Let $U\subset \nb_0$ be a closed disc containing $P_1,\ldots,P_r$ as interior points, and set $h:U\to \c$ the holomorphic  function  $h=\int_{P_1} \tau.$ By Theorem \ref{th:runge}, there exists $F\in \Fg_0(\nb_0)$ such that $|F|_U+h|<1$ and $(F|_U+h)_0 \geq  \prod_{j=1}^r P_j^{n_j+1}.$ This function solves the claim.
\end{proof}

Let us show that $\Gamma_1=\Gamma_2.$  Indeed, it is well known (see \cite{farkas}) that there exist $2 \nu+b-1$ cohomologically independent meromorphic 1-forms  in $\Wg_h(\nb_0)$  generating the first holomorphic De Rham cohomology group  ${\cal H}^1_{\text{hol}}(\nb_0)$  of $\nb_0.$ Recall that ${\cal H}^1_{\text{hol}}(\nb_0)$ is the quotient $\Wg_h(\nb_0)/\sim,$ where $\sim$ is the equivalence relation
$$\text{$\tau_1 \sim \tau_2$ if and only if $\tau_2-\tau_1=dh$ for some $h\in \Fg_h(\nb_0).$}$$  Thus, the map $H^1_{\text{hol}}(\nb_0) \to \c^{2 \nu+b-1}$,  $ [\tau] \mapsto \left( \int_{d} \tau \right)_{d  \in B_0},$ is a linear isomorphism. Assume that $\Gamma_1\neq \Gamma_2$ and take  $\tau \in \Wg_h(\nb_0)$ such that   $\int_{\Gamma_1-\Gamma_2} \tau \neq 0.$ By Claim \ref{ass:puntos},  we can find $F \in  {\Fg_h}(\nb_0)$ such that $(\tau+dF)_0 \geq  (d \hat{g})_0^3 (\hat{g})_\infty^2 \prod_{i=1}^a E_i^{2 m_i}.$ Set $h:=\frac{(\tau+dF) (\hat{g}^2-1)^2}{d \hat{g}}\in \Fg_h(S^c)$ and note that $(h) \geq (\hat{g}^2-1)_0^2 (d \hat{g})_0^2 \prod_{i=1}^a E_i^{2m_i}.$ By Theorem \ref{th:runge} and Remark \ref{re:runge},  $h$ lies in the closure of ${\cal L}_{1}$ in $(\Fg_h(S^c),\|\cdot\|_{0,S})$ and equation (\ref{eq:fun2}) gives that $ \int_{\Gamma_1-\Gamma_2} \tau+dF=\int_{\Gamma_1-\Gamma_2} \tau=0,$ a contradiction.

Coming back to equation (\ref{eq:fun1}) and using that $\Gamma_1=\Gamma_2,$ one has 
\begin{equation} \label{eq:fun3}
\int_{\Gamma_1} f (\frac{1}{\hat{g}^2}-1) d \hat{g}=0,\quad \mbox{for all} \; f \in {\cal L}_0.
\end{equation}
Let us see now that $\Gamma_1=0.$ Reason by contradiction and suppose that $\Gamma_1 \neq 0.$  As above, take  $\tau \in  \Wg_h(\nb_0)$ and $H \in   \Fg_h(\nb_0)$  such that   $\int_{\Gamma_1} \tau \neq 0$  and $(\tau+dH)_0 \geq  (d\hat{g})_0 (\hat{g}^2-1)_0 \prod_{i=1}^a E_i^{m_i}.$ Set  $t:=\frac{(\tau+dH) \hat{g}^2}{(\hat{g}^2-1) d \hat{g}}\in \Fg_h(S^c)$ and observe that $(t)\geq  (\hat{\phi}_3)_0^2\prod_{i=1}^a E_i^{m_i}.$  By Theorem \ref{th:runge} and Remark \ref{re:runge}, $t$ lies in the closure of ${\cal L}_0$ in $(\Fg_h(S^c),\|\cdot\|_{0,S})),$ hence from equation (\ref{eq:fun3}) we get that $\int_{\Gamma_1} (\tau+dH)=\int_{\Gamma_1} \tau=0,$ a contradiction.

Finally, equation (\ref{eq:fun}) and the fact that $\Gamma_1=\Gamma_2=0$ give that
\begin{equation} \label{eq:fun4}
\int_{\Gamma_3} f \hat{\phi}_3=0 \quad \mbox{for all}\; f \in {\cal L}.
\end{equation}
Reasoning as above, there exist $\tau \in  \Wg_h(N_0)$  and $G \in   \Fg_h(N_0)$ such that $\int_{\Gamma_3} \tau \neq 0$ and $(\tau+dG)_0 \geq  (\phi_3)_0.$  The function $v:=\frac{(\tau+dG)}{\hat{\phi}_3}$ lies in $\Fg_h^(S^c)$ and satisfies that $(v) \geq \prod_{i=1}^a E_i^{m_i}.$ By Theorem \ref{th:runge} and Remark \ref{re:runge}, $v$ lies in the closure of ${\cal L}$ in $(\Fg_h(S^c),\|\cdot\|_{0,S}))$  and equation (\ref{eq:fun4}) can be formally applied  to $v.$ We get that $\int_{\Gamma_3} \tau+dF=0,$ absurd. This contradiction proves the lemma.

\end{proof}


Now, we introduce  analytical deformation  and period operators for data $\{\eta_1^n,\eta_2^n\},$ $n\in \n.$  

\begin{definition}
For each $f \in {\cal L}$ and $n \in \n,$ set $\eta_j^n(f)=(1+j f)^2 \eta_j^n,$ $j=1,2,$ and  define $\phi_j^n(f),$ $j=1,2,3,$ like in equation \ref{eq:spinor}.
Set also $\Phi_n(f):=(\phi_j^n(f))_{j=1,2,3}.$
\end{definition}\label{def:defor}
It is clear that $\eta_j^n(f)-\eta_j^n,$ $\phi_k^n(f)-\phi_k^n \in \Wg_h(N_0),$  hence 
\begin{equation} \label{eq:rollo}
\Phi_n(f)-\hat{\Phi}\in \Wg_h(S^c)^3 \; \mbox{for all}\;  f \in {\cal L}\;\mbox{and}\;n \in \n.  
\end{equation}
Set ${\cal P}_n:{\cal L} \to \c^{\varsigma_0}, \quad {\cal P}_n(f)= \big(\int_{d} \Phi_n(f) -\hat{\Phi}\big)_{d  \in B_0},\quad n \in \n.$

Following Lemma \ref{lem:implicit}, let ${\cal U} \subset {\cal L}$ be a $\varsigma_0$-dimensional complex linear subspace such that $d {\cal P}_0({\cal U})=\c^{\varsigma_0},$ and fix a basis $\{f^j_d\,:\, d  \in B_0, \,j \in \{1,2,3\}\}$ of ${\cal U}.$ 

For the sake of simplicity, write ${{\bf f}_0}=[(f^j_d)_{d  \in B_0}]_{j=1,2,3} \in {\cal L}^{\varsigma_0}.$ 
For any  ${\bf x}=[(x^j_d)_{d  \in B_0}]_{j=1,2,3} \in \c^{\varsigma_0}$ and ${\bf h}=[(h^j_d)_{d  \in B_0}]_{j=1,2,3} \in {\cal L}^{\varsigma_0},$ write also ${\bf x} \cdot {\bf h}=\sum_{j=1}^3 [\sum_{d  \in B_0} x^j_d h^j_d].$ For each $n \in \n\cup \{0\}$ and ${\bf h} \in {\cal L}^{\varsigma_0},$  set ${\cal Q}_{n,{\bf h}}:\c^{\varsigma_0} \to \c^{\varsigma_0}$ for the vectorial degree two complex polynomial function given by
$${\cal Q}_{n,{\bf h}}({\bf x})={\cal P}_n( {\bf x} \cdot {\bf h}),$$ where we have made the convention ${\cal P}_0={\cal P}.$

By  Lemma \ref{lem:implicit}, ${\cal Q}_{0,{{\bf f}_0}}$ has non-zero Jacobian at the origin, hence we can find a closed Euclidean ball $K_0\subset \c^{\varsigma_0}$ centered at the origin such that ${\cal Q}_{0,{{\bf f}_0}}|_{K_0}: K_0 \to {\cal Q}_{0,{{\bf f}_0}}(K_0)$ is a biholomorphism. Moreover, since ${\cal Q}_{0,{{\bf f}_0}}({\bf 0})=0$ then ${\cal Q}_{0,{{\bf f}_0}}(K_0)$ contains the origin as an interior point. Since $\{{\cal Q}_{n,{{\bf f}_0}}\}_{n\in \n}\to {\cal Q}_{0,{{\bf f}_0}}$ uniformly on compact subsets of $\c^{\varsigma_0}$ and the convergence is analytical, then ${\cal Q}_{n,{{\bf f}_0}}|_{K_0}: K_0 \to {\cal Q}_{n,{{\bf f}_0}}(K_0)$ is a biholomorphism, and ${\cal Q}_{n,{{\bf f}_0}}(K_0)$ is an Euclidean ball containing the origin is an interior point as well, $n$ large enough (without loss of generality, for all $n$). Let ${\bf x}_n\in K_0$ denote the unique point such that ${\cal Q}_{n,{{\bf f}_0}}({\bf x}_n)=0,$ and set  ${h}_n:={\bf x}_n \cdot {\bf f}_0\in {\cal L},$ $n\in \n.$

The sequence $\{h_n\}_{n\in \n}$ solves the second phase of the program.  Indeed, one has that $\{\Phi_n({h}_n)\}_{n\in \n}$ uniformly approximates $\hat{\Phi}$ on $S,$ $\{\Phi_n({ h}_n)\}_{n\in \n}$ have no real periods on $S,$ and $\{\Phi_n({ h}_n)\}_{n\in \n}$ induce the same flux map as $\hat{\Phi}$ (just notice that $\hat{\Phi}-\Phi_n({h}_n)$ is exact on $S$ for all $n$). The second tentative of solution for the Fundamental Approximation Theorem is to define $Y_n:\nb \to \r^3,$ $Y_n:=\mbox{Re} \left(\int\Phi_n(h_n)\right)$ for all $n\in \n.$ However, the 1-forms $\eta_n^1(h_n)$ and  $\eta_n^2(h_n)$ could have common zeros in $\nb,$ and consequently $Y_n$ could fail to be an immersion. Even more, we have no control over the behavior of $Y_n$ on the punctures of $\nb_0.$  

To overcome this difficulties, we have to devise a more sophisticated deformation procedure. This is the content of the following phase.

\subsubsection{Third phase: proving the theorem.} \label{subsub:phase3}
Let us keep the notations of the previous paragraphs. 

To finish the proof, we are going to reproduce the previous program but replacing  ${\bf f}_0$ for a suitable basis ${\bf f}_n$ of ${\cal U}$ depending on $n\in \n$ (see Lemma \ref{lem:4} below).

Up to choosing a smaller ball $K_0\subset \c^\varsigma_0,$  in the sequel we will assume that 
\begin{equation} \label{eq:ajuste}
\|{\bf x} \cdot {\bf f}_0\|_{0,S}<1 \; \mbox{for all}\; {\bf x}\in K_0.
\end{equation}

\begin{lemma} \label{lem:4}
We can find $\{{\bf f}_n\}_{n \in \n} \subset {\cal L}^{\varsigma_0}$ such that:
\begin{enumerate}[(i)]
  \item $\{{\bf f}_n|_{S^c}\}_{n \in \n} \to {\bf f}_0|_{S^c}$ in the $\mathcal{C}^0(S)$-topology.
  \item $\eta_j^n({\bf x} \cdot {\bf f}_n)$ has a pole at $E_k$ for all $k\in \{1,\ldots,a+b\},$ ${\bf x}  \in \c^{\varsigma_0}$ and $n \in \n,$ $j=1,2,$
  \item  $\sum_{j=1}^2 |\eta_j^n({\bf x} \cdot {\bf f}_n)|$ never vanishes on $\nb$ for all $n \in \n$ and ${\bf x}\in K_0.$
\end{enumerate}
As a consequence, $\{\eta_j^n({\bf x} \cdot {\bf f}_n)|_{S}\}_{n \in \n} \to \hat{\eta}_j({\bf x} \cdot {\bf f}_0)$ in the $\mathcal{C}^0(S)$-topology and  $\{{\cal Q}_{n,{\bf f}_n}\}_{n \in \n} \to {\cal Q}_{0,{{\bf f}_0}}$ uniformly on $K_0.$
\end{lemma}
\begin{proof} By definition, it is clear that $\eta_j^n(f)$ has poles at  $E_k$ for all $k\in\{1,\ldots,a\}$ and $f \in {\cal L}.$ By Lemma \ref{lem:sp},   $\{\eta_j^n|_{S^c}-\hat{\eta}_j\}_{n \in \n} \to 0$ in the $\mathcal{C}^0(S)$-topology,  $(\eta_j^n|_{S^c})=(\hat{\eta}_j),$ and $\eta_j^n$ never vanishes on $\partial(S),$ $j=1,2,$ for all $n.$ Let $C_n$ be a finite collection of closed discs in $N_0-S^c$ containing all the zeros of $\eta_1^n$ and $\eta_2^n$ in $N_0-S^c$ and meeting all the bounded components of $N^c-S^c.$ Obviously, $S^c \cup C_n$ is admissible in the open Riemann surface $\nb_0.$

For each $d  \in B_0,$  $n \in \n,$ and $j \in \{1,2,3\},$ set $\hat{f}_d^{j,n}:S^c  \cup C_n \to \c,$ $\hat{f}_d^{j,n}|_{S^c}=f^j_d,$ $\hat{f}^{j,n}_d|_{C_n}=0.$

By Theorem \ref{th:runge} and similar arguments to those used in the proof of Lemma \ref{lem:sp}, we can find  a sequence $\{f^{j,n}_{d}(m)\}_{m \in \n}$ in ${\cal L}$ satisfying that
\begin{itemize}
\item  $\{f^{j,n}_{d}(m)|_{S^c\cup C_n}\}_{m \in \n}\to \hat{f}^{j,n}_d$ in the $\mathcal{C}^0(S^c \cup C_n)$-topology, and
\item  the sequence of pole multiplicities $\{\mbox{Ord}_{E_{a+k}} \big(f^{j,n}_{d}(m)\big)\}_{m \in \n}$ is divergent for all $k \in \{1,\ldots,b\}.$
\end{itemize}
Up to  subsequences, we can assume that:
\begin{equation}\label{eq:polos}
\mbox{Ord}_{E_{a+k}}\big( f^{j_1,n}_{d_1}(m) \big)\neq \mbox{Ord}_{E_{a+k}}\big( f^{j_2,n}_{d_2}(m)\big)\; \mbox{provided that}\; (d_1,j_1) \neq (d_2,j_2).
\end{equation}

Set  ${\bf f}_n(m)=[(f^{j,n}_d(m))_{d  \in B_0}]_{j=1,2,3},$ $m \in \n,$ and take a divergent sequence $\{m_n\}_{n \in \n}\subset \n$ such that  $\{{\bf f}_n(m_n)|_{S^c}\}_{n \in \n} \to {\bf f}_0|_{S^c}$ in the $\mathcal{C}^0(S^c)$-topology and $\big\{\max_{C_n}|f^{j,n}_d(m_n)| \big\}_{n \to \n} \to 0$ for all $d\in B_0$ and $j\in \{1,2,3\}.$

Set ${\bf f}_n={\bf f}_n(m_n)$ for all $n \in \n,$ and let us show that $\{{\bf f}_n\}_{n \in \n}$ solves the claim. 

Items $(i)$ is  obvious, and item $(ii)$ follows from equation (\ref{eq:polos}) and the facts that $\eta_j^n$ has a pole at $E_k$ for all $k \in \{a+1,\ldots,a+b\},$ $\eta_j^n(f)$ has a pole at  $E_j$ for all $j\in\{1,\ldots,a\},$ and $f \in {\cal L}.$  

Let us check $(iii).$ Taking into account (i) and equation (\ref{eq:ajuste}), and removing finitely many terms of $\{{\bf f}_n\}_{n \in \n}$ if necessary, we can assume that:
\begin{enumerate}[(a)]
  \item $\|{\bf x} \cdot {\bf f}_n\|_{0,S}<1$  for all ${\bf x}\in K_0,$  and so $1+{\bf x} \cdot {\bf f}_n$ and $2+{\bf x} \cdot {\bf f}_n$ never vanish on $S^c$ for all ${\bf x}\in K_0,$ 
  \item $1+{\bf x} \cdot {\bf f}_n$ and $2+{\bf x} \cdot {\bf f}_n$ never vanish on $C_n$ for all $n$ and ${\bf x}\in K_0.$
\end{enumerate}
Since $|\eta_1^n|+|\eta_2^n|$ never vanishes on $\nb,$ $(a)$ and $(b)$ give that $\sum_{j=1}^2 |\eta_j^n({\bf x} \cdot {\bf f}_n)|$ never vanishes on  $S \cup C_n$ for all $n.$  
Taking into account that  $|1+f|+|1+2 f|$ never vanish on $\nb$ for all $(n,f)\in  \n\times {\cal L},$ and the choice of $C_n,$ we deduce that   $\sum_{j=1}^2 |\eta_j^n({\bf x} \cdot {\bf f}_n)|$ never vanishes on $\nb-(S \cup C_n)$ as well, and  we are done.
\end{proof}

With the help of this lemma, we can tackle the decisive part of the proof.\\

For the sake of simplicity, write ${\cal Q}_n= {\cal Q}_{n,{\bf f}_n},$ $n \in \n\cup \{0\}.$ At this point, we reproduce the previous program once again. Since the coefficients of the vectorial polynomial functions $\{{\cal Q}_n\}_{n \in \n}$  converge to the ones of ${\cal Q}_0$ (take into account Lemma \ref{lem:4}),   ${\cal Q}_n|_{K_0}: K_0 \to {\cal Q}_n(K_0)$ is a biholomorphism and ${\cal Q}_n(K_0)$  contains the origin as an interior point, $n$ is large enough (up to removing finitely many terms, for all $n$). 

Let ${\bf y}_n \in K_0$ denote the unique point satisfying ${\cal Q}_n({\bf y}_n)={\bf 0},$ and notice that $\lim_{n \to \infty}{\bf y}_n ={\bf 0}.$  Set  $\rho_j^n:=\eta_j^n({\bf y}_n \cdot {\bf f}_n ),$ $j=1,2,$ $\psi_k^n:=\phi_k^n({\bf y}_n \cdot {\bf f}_n ),$ $k=1,2,3,$ and define $$Y_n:N \to \r^3,\quad Y_n(P)=X(P_0)+\mbox{Re} \int_{P_0}^P (\psi_k^n)_{k=1,2,3},\; n \in \n,$$ where $P_0$ is any point of $S.$

Now we can prove that $\{Y_n\}_{n\in \n}$ is the solution for the first part of the Fundamental Approximation Theorem. 

Indeed, by (\ref{eq:rollo}) and the choice of ${\bf y}_n$, $\psi_k^n-\hat{\phi}_k$ is an exact 1-form in $\Wg_h(S^c)$ and  $Y_n$ is well defined. Moreover,  Lemma \ref{lem:4} and  Osserman's theorem imply that $Y_n \in  \mathcal{M}_{X_\sigma}(\nb).$  As  $\rho_j^n|_{S^c}= \big(1+j ({\bf y}_n \cdot {\bf f}_n|_{S^c}) \big)^2 (\eta_j^n|_{S^c}-\hat{\eta}_j)+ \big(1+j ({\bf y}_n \cdot {\bf f}_n)|_{S^c} \big)^2 \hat{\eta}_j,$  then Lemma \ref{lem:sp},  Lemma \ref{lem:4} and the fact that $\{{\bf y}_n\}_{n \in \n} \to 0$ give that $\{\rho_j^n|_S-\hat{\eta}_j\}_{n \to \n} \to 0$ in the $\mathcal{C}^0(S)$-topology. Therefore $\{\Rg_S(Y_n)\}_{n\in \n} \to X_\sigma$ in the ${\cal C}^1(S)$-topology, proving the first part of the Theorem \ref{th:density}.\\

The second and final part of the theorem is a direct application of the previous ideas and the well-known Jorge-Xavier theorem (see  \cite{jorge-xavier}). Let  $V$ be  a closed tubular  neighborhood of $S.$   Let $L_n$ be a Jorge-Xavier type labyrinth in $V^\circ-S$  adapted to $C$ and $ \psi_3^n,$ that is to say, a finite collection of pairwise disjoint closed discs in  $V-S$ such that $\int_{\gamma} | \psi_3^n| > C$ for any compact arc $\gamma \subset  V-L_n$ connecting $\partial(S)$ and $\partial(V)$ (see \cite{jorge-xavier} or \cite{nadi}). Consider another Jorge-Xavier type labyrinth $L'_n$   obtained as a small closed tubular  neighborhood of $L_n$ in $V^\circ-S.$  By Theorem \ref{th:runge}, there is  $\{h_{n,m}\}_{m \in \n} \subset \Fg_h(\nb_0)$ such that $|h_{m,n}|<1/m$ on $S^c,$ $|h_{m,n}-m|<1/m$ on $L'_n$ and $(h_{n,m})_0 \geq \prod_{j=1}^a E_j^{m_i},$ $m \in \n.$ 

Consider on $\nb$ the spinorial data $\varphi_1^{n,m}=e^{-h_{n,m}}  \eta_1^n,$ $\varphi_2^{n,m}=e^{h_{n,m}}  \eta_2^n,$   and  their associated Weierstrass data
$$\text{$\tau_1^{n,m}=1/2 (\varphi_1^{n,m}-\varphi_2^{n,m}),$ $\tau_2^{n,m}=i/2 (\varphi_1^{n,m}+\varphi_2^{n,m})$ and $\tau_3^{n,m}=\psi_3^{n}.$}$$
For any $f \in {\cal L}$ put $\varphi_1^{n,m}(f)=e^{-h_{n,m}}  \eta_1^n(f),$ $\varphi_2^{n,m}(f)=e^{h_{n,m}}  \eta_2^n(f),$ and call 
$$\text{$\tau_1^{n,m}(f)=1/2 (\varphi_1^{n,m}(f)-\varphi_2^{n,m}(f)),$ $\tau_2^{n,m}=i/2 (\varphi_1^{n,m}(f)+\varphi_2^{n,m}(f))$ and $\tau_3^{n,m}(f)=\psi_3^{n}(f).$}$$
Define the period operator ${\cal Q}_{n,m}:\c^{\varsigma_0} \to \c^{\varsigma_0},$ ${\cal Q}_{n,m}({\bf x})=\big[(\int_{d} \tau_j^{n,m}({\bf x}\cdot {\bf f}_n  )-\hat{\phi}_j)_{d  \in B_0}\big]_{j=1,2,3}.$ One has  $\{{\cal Q}_{n,m}\}_{m \in \n} \to {\cal Q}_n$  uniformly on compact subsets of $\c^{\varsigma_0},$  ${\cal Q}_{n,m}|_{K_0}: K_0 \to {\cal Q}_{n,m}(K_0)$ is a biholomorphism, and  $0 \in {\cal Q}_{n,m}(K_0)-\partial \big({\cal Q}_{n,m}(K_0) \big)$ for large enough  $m$ (without loss of generality  for all $m$). Therefore,  $\lim_{m \to \infty}{\bf y}_{n,m} ={\bf y}_n,$ where ${\bf y}_{n,m} \in K_0$ is the unique point satisfying ${\cal Q}_{n,m}({\bf y}_{n,m})={\bf 0}.$

Call $\psi_j^{n,m}=\tau_j^{n,m}({\bf y}_{n,m} \cdot {\bf f}_n),$ $j=1,2,3,$ fix $P_0\in S$ and set $$Y_{n,m}:\nb \to \r^3,\quad Y_{n,m}(P)=X(P_0)+\mbox{Re} \int_{P_0}^P (\psi_k^{n,m})_{k=1,2,3}.$$ Note that $Y_{n,m}$ is  well defined, has no branch points (take into account Lemma \ref{lem:4}), $\Rg_S(Y_{n,m})\in \mathcal{M}^*(S),$   and $\Rg_S(Y_{n,m})$  and $X_\sigma$ are  flux equivalent on $S,$ $m \in \n.$ Moreover,  $\{\Rg_S(Y_{n,m})\}_{m \in \n} \to \Rg_S(Y_n)$ in the ${\cal C}^1(S)$-topology for all $n.$

From the choice of $L_n'$ and the fact  $\{e^{h_{m,n}}\}_{m \in \n}\to \infty$ uniformly on $L'_n,$ one has that $d_{Y_{n,m}}(S,\partial({V}))> C$ for large enough $m$ (depending on $n$), where $d_{Y_{n,m}}$ is the intrinsic distance in $\nb$ associated to $Y_{n,m}.$  Since $\{\Rg_S(Y_n)\}_{n \in \n} \to X_\sigma$ in the ${\cal C}^0(S)$-topology,  for each $n$ we can find $m_n\in \n$  such that the immersions $H_n=Y_{n,m_n},$ $n\in \n,$ satisfy:
\begin{itemize}
\item  $d_{H_{n}}(S,\partial({V}))> C,$ where $d_{H_n}$ is the intrinsic distance in $\nb$ associated to $H_n,$  and
\item  $\{\Rg_S(H_n)\}_{n \in \n} \to X_\sigma$ in the ${\cal C}^1(S)$-topology.
\end{itemize}

Unfortunately, $H_n$ is not necessarily of FTC, and we have to work a little more.  Applying the first part of the theorem to $H_n|_V$ (notice that $V$ is admissible in $\nb$), there exists $\{Z_{n,j}\}_{j \in \n} \subset \mathcal{M}(\nb)$ such that $\{Z_{n,j|_V}\}_{j \in \n}\to  H_n|_V$ in the ${\cal C}^0(V)$-topology, and $\Rg_S(Z_{n,j})$ and $\Rg_S(H_n)$ are flux equivalent  for all $j.$ In particular, $\{Z_{n,j}\}_{j \in \n} \subset \mathcal{M}_{X_\sigma}(N)$ for all $n\in \n.$   Furthermore, without loss of generality we can also suppose that $d_{Z_{n,j}}(S,\partial({V}))> C$ for all $j$ and $n.$

Since  $\{\Rg(H_n)\}_{n \in \n} \to X_\sigma$ in the $\mathcal{C}^1(S)$-topology,  a standard diagonal process provides a sequence $\{Z_n\}_{n \in \n} \subset \{Z_{n,j}\,|\; n,\,j\in \n\} \subset   \mathcal{M}_{X_\sigma}(N)$ such that $\{\Rg_S(Z_n)\}_{n \in \n} \to X_\sigma$ in the ${\cal C}^1(S)$-topology and  $d_{Z_{n}}(S,\partial({V}))\geq C$ for all $n,$ concluding the proof.

\subsection{General version of the Fundamental Approximation Theorem.} \label{subsec:parabo}

In this subsection we obtain the general version of the Fundamental Approximation Theorem for arbitrary open Riemann surfaces and admissible subsets.

Let us start with the following

\begin{lemma} \label{lem:opera}
Let $\nb$ be an open Riemann surface, let $S=M\cup\beta$ be a (possibly non-connected) admissible subset in $\nb,$ and let $V$ be an admissible region in $\nb$ of finite conformal type containing $S.$ 
Let $X_\sigma$ be a marked immersion in $\mathcal{M}^*(S),$ let   $q:\mathcal{H}_1(V,\z) \to \r^3$ be  a  group homomorphism satisfying that $q|_{\mathcal{H}_1(S,\z)}=p_{X_\sigma},$ and fix arbitrary constants $C>0,$ $\epsilon>0.$ 

Then there exists $Y \in \mathcal{M}(V)$ such that $\|Y-X_\sigma\|_{1,S}\leq \epsilon,$ $d_{X}(S,\partial (V)) \geq C,$ and  $p_Y=q.$  
\end{lemma}
\begin{proof}  By basic topology, we can find a finite collection $\gamma\subset V$ of Jordan arcs such that $S_0=S\cup \gamma$ is a connected admissible subset in $\nb$ and $j_*:\mathcal{H}_1(S_0,\z)\to \mathcal{H}_1(V,\z)$ is an isomorphism, where $j:S\to V$ is the inclusion map. This simply means that $V-S_0$ consists of a finite collection of once punctured discs and conformal annuli.  

If $V\neq \nb,$ consider  a closed tubular  neighborhood $V_0$ of $V$ in $\nb$ and  a conformal compactification $R$  of $V_0.$  Recall that $R-V_0$ consists of a finite family $U_1,\ldots,U_r$ of pairwise disjoint open discs. Moreover, if we fix $P_j\in U_j$ for each $j,$  $S_0$ becomes  an admissible subset of $R_0:=R-\{P_1,\ldots,P_r\},$ and $R_0-S_0$ consists of $r$ pairwise disjoint once punctured open discs. If $V=\nb,$ simply set $V_0=R_0=V$ and $R=\nb^c.$

Construct ${X_0}_{{\sigma_0}}\in \mathcal{M}^*(S_0)$  satisfying that $X_0|_{S}=X$ and $\sigma_0|_\beta=\sigma.$ By Theorem \ref{th:density}, there exists $Z\in \mathcal{M}(R_0)$ such that $\|Z-X_0\|_{1,S_0}<\epsilon$ and $p_Z|_{\mathcal{H}_1(S_0,\z)}=p_{{X_0}_{\sigma_0}}.$ The immersion $Y:=Z|_{V}$ solves the lemma.

\end{proof}

\begin{theorem}[General Approximation Theorem]  \label{th:parabo}
Let $\nb$ be an open Riemann surface, and let $S$ be a possibly non connected  admissible subset in $\nb.$  Let $X_\sigma \in \mathcal{M}^*(S)$ and let  
 $q:\mathcal{H}_1(\nb,\z) \to \r^3$  be a group morphism such that $q|_{\mathcal{H}_1(S,\z)}=p_{X_\sigma}.$

Then there exists a  sequence  $\{Y_n\}_{n \in \n}\in \mathcal{M}(\nb)$ such that  $\{\Rg_S(Y_n)\}_{n \in \n} \to X_\sigma$ in the ${\cal C}^1(S)$-topology and $p_{Y_n}=q$ for all $n.$
\end{theorem}
\begin{proof}  It suffices to prove that for any $\epsilon>0$ there is $Y \in \mathcal{M}(\nb)$ such that $\|Y-X_\sigma\|_{1,S}\leq \epsilon$ and $p_Y=q.$ 

If $\nb$ is of finite conformal type,   the theorem follows from Lemma \ref{lem:opera}.  In the sequel we will suppose that $\nb$ is not  of finite conformal type, or equivalently that  $\nb^c$ is non-compact. Write $E=\nb^c-\nb.$

Consider an  exhaustion $\hat{N}_1 \subset  \hat{N}_2 \subset  ...$ of $\nb^c$ by  compact regions such that  
\begin{itemize}
\item $\hat{N}_0:=S^c\subset  \hat{N}_1^\circ,$ 
\item $ \hat{N}_j$ is admissible in the open Riemann surface $\nb^c$ for all $j\geq 1,$ and
\item $ \hat{N}_j \subset  \hat{N}_{j+1}^\circ$ and $E_j:=E\cap  \hat{N}_j\subset  \hat{N}_j^\circ$ for all $j\geq 1.$
\end{itemize}
Call $N_j=-\hat{N}_j-E_j,$ $j\geq 1,$ and set $Y_0=X_\sigma.$

Using Lemma \ref{lem:opera} in a recursive way, one  can  construct $Y_j \in \mathcal{M}({N}_j),$ $j \geq 1,$ satisfying that:
\begin{enumerate}[(i)]
\item  $\|Y_{j+1}-Y_{j}\|_{1,N_j} \leq \epsilon/2^{j+1}$ and  $p_{Y_j}=q|_{\mathcal{H}_1(N_j,\r)},$  for all $j \geq 0.$  
\item  $d_{Y_{j+1}}\big(Y_{j+1}(N_{j}), Y_{j+1} \big(\partial(N_{j+1})\big)\big) \geq 1,$  where  $d_{Y_{j+1}}$ means intrinsic distance in $\nb$ with respect to $Y_{j+1},$  $j\geq 0.$
\end{enumerate}
Let $Y:\nb\to\r^3$ be the possibly branched minimal immersion given by $Y|_{N_j}=\lim_{m \to \infty} Y_m|_{N_j},$ $j \in \n,$  and note that $\lim_{m \to \infty} \|Y_m- Y\|_{1,N_j}=0$ for all $j$ and  $\|Y-X_\sigma\|_{1,S}\leq \epsilon.$  

Let us show that $Y$ has no branch points. 

Without loss of generality, we will suppose that $X$ is non-flat on the regions of $S$ (use similar ideas to those in  Proposition \ref{pro:nonflat}).  Up to choosing $\epsilon$ small enough, the inequality $\|Y-X_\sigma\|_{1,S}\leq \epsilon$ implies that $Y$ is non-flat as well. Let $(g_m,\phi_3^m)$ denote the Weierstrass data of $Y_m,$ $m \in \n,$ and likewise call $(g,\phi_3)$ the ones of $Y.$ Obviously, $\{g_m,\phi_3^m)\}_{m \in \n} \to (g,\phi_3)$ uniformly on compact subsets of $\nb^c.$  Take an arbitrary $P_0 \in \nb,$ and consider $j_0 \in \n$ such that  $P_0 \in N_{j_0}^\circ.$ Up to a rigid motion, $g(P_0) \neq 0,$ $\infty,$ hence we can find a closed disc $D \subset N_{j_0}$ such that $P_0 \in D^\circ$ and $g_m|_D,$ $m \in \n,$ $g|_D$ are holomorphic and never vanishing. Since $Y_m$ has no branch points, $\phi_3^m$ has no zeros on $D$ for all $m.$  By Hurwith theorem, either $\phi_3=0$ of $\phi_3$ has no zeros on $D$ as well. In the first case the identity principle would give $\phi_3=0$ on $\nb,$ contradicting that $Y$ is non-flat. Therefore, $\phi_3$ has no zeros on $D$ and $Y|_D$ has no branch points. Since $P_0$ is an arbitrary point of $\nb,$ $Y$ is a conformal minimal immersion.

Finally, let us see that $Y$ is complete and of WFTC. By Osserman's theorem, the Gauss map of $Y_j$ extends meromorphically to  $\hat{N}_j,$  $j \in \n.$  Since $\|Y_j-Y\|_{1,N_j}$ is finite, then   Weierstrass data of $Y$ extends meromorphically to $\nb^c$ as well and $Y|_{N_j}$ is complete and of  finite total curvature  for any $j.$ It remains to check that $Y$ is complete. Indeed, obviously those curves in $\nb$ diverging to a puncture in $E$ have infinite intrinsic length with respect to $Y.$ By  item $(ii),$  any curve in $\nb$ diverging in $\nb^c$ has also infinite intrinsic length. This shows that $Y$ is complete and lies in $\mathcal{M}(\nb).$ Since $p_{Y}=q,$ this completes the proof.   
\end{proof}
For any $X \in  \mathcal{M}(\nb)$ with $p_X=0$ and $\theta \in \partial(\d),$ we set $X_\theta=\mbox{Re}\big( \int  \theta \cdot \partial_z X \big)$ and call  $\{X_\theta \;:\; \theta \in \partial(\d)\} \subset  \mathcal{M}(\nb)$ as the family of {\em associated minimal immersions} of $X.$ The next corollary generalizes Pirola's results in \cite{pirola}:
\begin{corollary} \label{co:arbitrary}
For any open Riemann surface $\nb,$ there exists $Y \in  \mathcal{M}(\nb)$ such that all its associated immersions are well defined. In particular, the space $ \mathcal{M}(\nb)\neq \emptyset.$
\end{corollary}
\begin{proof}
Fix a closed disc $D \subset \nb$ and an immersion $X\in  \mathcal{M}(D).$ By Theorem \ref{th:parabo}, there is $\{Y_n\}_{n \in \n} \subset  \mathcal{M}(\nb)$ such that $\{Y_n|_D\}_{n \in \n} \to X$ in the ${\cal C}^0(D)$-topology and $p_{Y_n}=0.$   The corollary follows straightforwardly.
\end{proof}

{\bf FRANCISCO J. LOPEZ} \newline
Departamento de Geometr\'{\i}a y Topolog\'{\i}a \newline
Facultad de Ciencias, Universidad de Granada \newline
18071 - GRANADA (SPAIN) \newline
e-mail: fjlopez@ugr.es

\end{document}